\documentclass{gtart_a}
\pdfoutput=1
\usepackage{pinlabel}

\title[Stabilization in the braid groups II]{Stabilization in the braid groups II:\\Transversal simplicity of knots}

\author{Joan S Birman}
\givenname{Joan S}
\surname{Birman}
\address{Department of Mathematics\\Barnard College\\Columbia University\\\newline
2990 Broadway\\New York, NY 10027\\USA}
\email{jb@math.columbia.edu}
\urladdr{}

\author{William W Menasco}
\givenname{William W}
\surname{Menasco}
\address{Department of Mathematics\\University at Buffalo\\
Buffalo, NY  14260\\USA}
\email{menasco@math.buffalo.edu}
\urladdr{}

\volumenumber{10}
\issuenumber{}
\publicationyear{2006}
\papernumber{36}
\lognumber{0629}
\startpage{1425}
\endpage{1452}

\doi{}
\MR{}
\Zbl{}

\keyword{contact structures}
\keyword{tight}
\keyword{transversal knot type}
\keyword{3-braids}
\keyword{flypes}
\keyword{Bennequin invariant}
\keyword{transversally simple}
\subject{primary}{msc2000}{57M25}
\subject{secondary}{msc2000}{57M50}
\subject{secondary}{msc2000}{53C15}

\received{24 June 2005}
\revised{}
\accepted{28 June 2006}
\published{4 October 2006}
\publishedonline{4 October 2006}
\proposed{Rob Kirby}
\seconded{David Gabai, Cameron Gordon}
\corresponding{}
\editor{}
\version{}

\arxivreference{math.GT/0310280}

\AtBeginDocument{\let\bar\wbar\let\hat\what}
\makeop{top}
\makeop{min}

\newtheorem{theorem}{Theorem}
\newtheorem{proposition}{Proposition}
\newtheorem{lemma}{Lemma}

\theoremstyle{remark}

\def\reals{\mathbb R}

\let\pf\proof
\let\endpf\endproof
\def\be{\begin{enumerate}}
\def\ee{\end{enumerate}}
\def\bi{\begin{itemize}}
\def\ei{\end{itemize}}
\def\ms{\medskip}
\def\fib{{\bf H}}
\def\axis{{\bf A}}

\def\cA{\cal A}

\def\cB{{\cal B}}
\def\cC{{\cal C}}
\def\cD{{\cal D}}
\def\cE{{\cal E}}

\def\cG{{\cal G}}

\def\cTX{{\cal TX}}
\def\cT{{\cal T}}

\def\cX{{\cal X}}

\def\e{\epsilon}
\def\ca{{\kern.35em{\cal C}\kern-.35em{\cal A}}}
\def\pa{{\kern.35em{\cal P}\kern-.35em{\cal A}}}
\def\ta{{\kern.35em{\cal T}\kern-.35em{\cal A}}}
\def\td{{\kern.25em{T}\kern-.25em{\Delta}}}
\def\glow{\gamma_-}
\def\ghigh{\gamma_+}

\def\alow{a_-}
\def\ahigh{a_+}
\def\d{\delta}

\def\g{\gamma}
\def\Dhigh{\cD_+}
\def\Dlow{\cD_-}
\def\DDhigh{\Delta_+}
\def\DDlow{\Delta_-}
\def\Xhigh{X_+}
\def\Xlow{X_-}

\def\TXhigh{TX_+}
\def\TXlow{TX_-}
\def\TXmed{TX_0}

\def\bt{{\bf t}}

\def\e{{\epsilon}}
\def\TD{T\Delta}

\begin{document}

\begin{abstract}
The main result of this paper is a negative answer to the question:
are all transversal knot types transversally simple?  An explicit
infinite family of examples is given of closed 3--braids that define
transversal knot types that are not transversally simple. The method
of proof is topological and indirect.
\end{abstract}

\maketitle

\section{Introduction and description of results}
\label{section:introduction}

This paper is about knots which are transverse to the standard tight
contact structure in $\reals^3$, and about the ways in which
topological information about braids can be used to learn new things
about these `transversal knots'.  Our approach to contact topology had
its origins in the foundational paper of D Bennequin
\cite{Bennequin}. The results in this paper are direct outgrowths of
our work on the {\em Markov Theorem without Stabilization} (MTWS)
\cite{BM-MTWS}. 

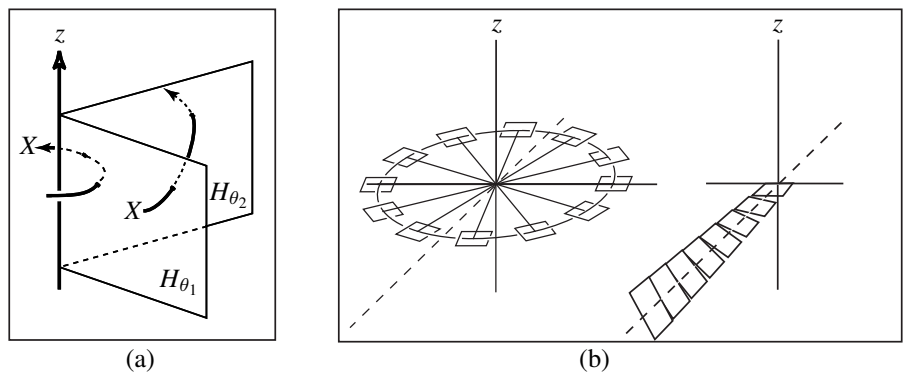
\begin{figure}[ht!]
\labellist\small
\pinlabel $z$ [b] at 66 690
\pinlabel $z$ [b] at 321 696
\pinlabel $z$ [b] at 485 696
\pinlabel (a) [t] at 113 521
\pinlabel (b) [t] at 378 521
\pinlabel* $X$ [r] at 53 635
\pinlabel* $X$ [r] at 114 597
\pinlabel $H_{\theta_1}$ at 136 555
\pinlabel $H_{\theta_2}$ at 165 606
\endlabellist
\centerline{\includegraphics[width=.9\hsize]{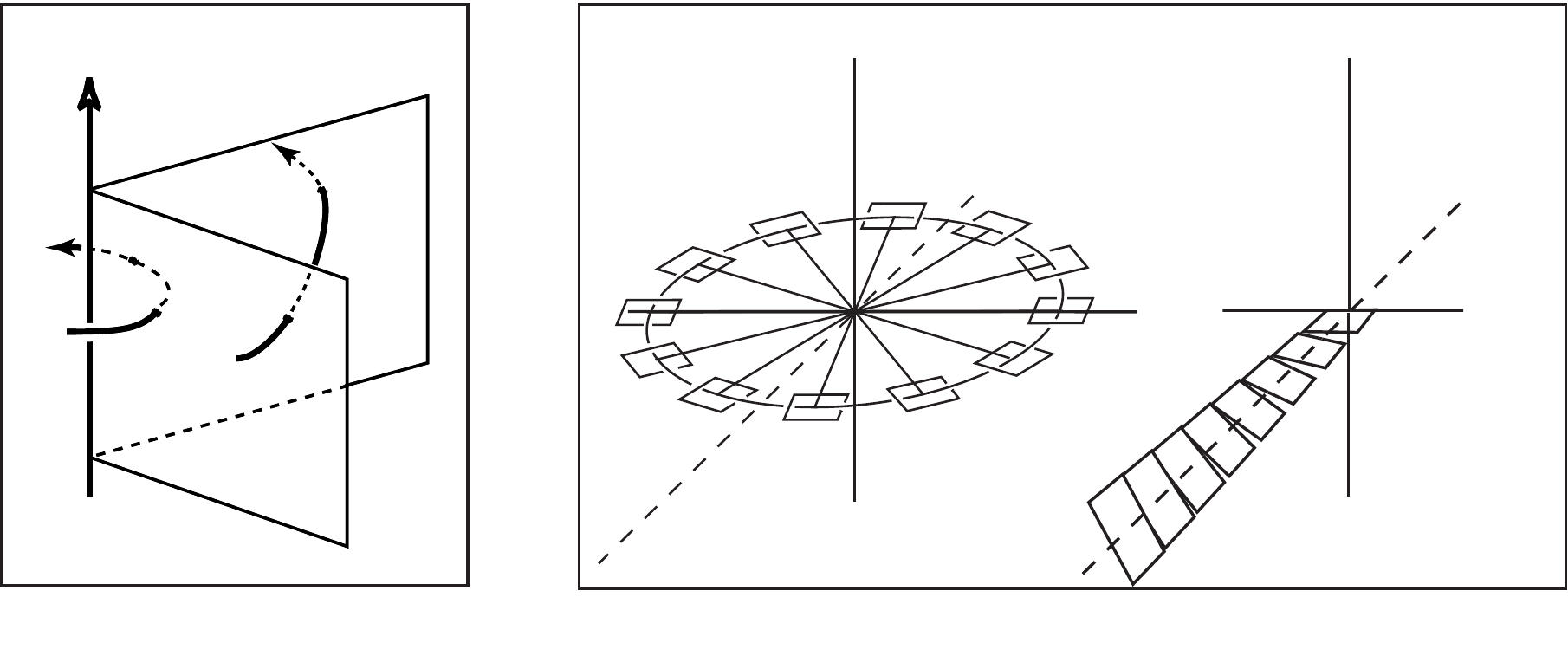}}  
\caption{(a)\qua Braid structure on $\reals^3$;\qua (b)\qua 2--planes in the standard contact structure on $\reals^3$.
Note that these 2--planes are invariant under rotation about $z$--axis and translation along
the  $z$--axis (which we use as the braid axis $\axis$).}
\label{figure:braid-contact}
\end{figure}

A {\em knot} in oriented $\reals^3$
is the image $X$ of an oriented circle $S^1$ under a smooth embedding $e\co S^1 \to \reals^3$.
Viewing $S^3$ as $\reals^3 \cup {\infty}$, we also can think of $X$ as being a knot in $S^3$.
The  {\em topological knot type} $\cX$ of $X$ is its equivalence
class under smooth isotopy of the pair $(X,S^3)$. A representative $X \in \cX$ is said to be a
{\em closed braid}
if there is an unknotted curve $\axis \subset (S^3 \setminus X)$, the {\em axis}, 
and a {\em choice of fibration $\fib$} of the open solid torus $S^3 \setminus \axis$ by meridian discs
$ \{ H_{\theta}; \ \theta \in [0,2\pi) \}$ which we call {\em fibers},  such that $X$
intersects every disc fiber  $H_\theta$ transversally. We call the pair
$(\fib,\axis)$ a {\em braid structure} on $\reals^3$. If we regard
$\axis$ as the
$z$--axis in
$\reals^3$ and the fibers of
$\fib$ as half-planes through $\axis$, there are natural cylindrical coordinates $(\rho,\theta,z),
\
\theta \in [0,2\pi)$ in
$\reals^3$. Parameterizing $X$ by $(\rho (t), \theta (t), z(t)), \ t \in [0,2\pi)$, $X$ will be a closed
braid if $\rho(t) > 0$ and 
$d\theta/dt >0$ at every point $(\rho (t), \theta (t), z(t))\in X$.  See \fullref{figure:braid-contact}(a), where two half-planes and two fragments of a closed braid $X$ are
depicted.  A {\em braid isotopy} of a closed braid is an isotopy that is transverse to the braid
fibration and in the complement of the axis.  The {\em braid isotopy type} of the closed braid $X$ is its 
equivalence class under braid isotopy.

Using the same cylindrical coordinates, the standard tight  {\em contact structure} $\xi$ in oriented
$\reals^3$  is the kernel of the 1--form $\alpha = \rho^2 d \theta + d z$. Thinking of $\xi$ as attaching a
2--plane to each point in $\reals^3$, \fullref{figure:braid-contact}(b) shows some of the 2--planes associated to $\alpha$, for comparison of the braid
structure and the contact structure. The field of 2--planes so-obtained is nowhere integrable. 

An oriented knot
$X$ in contact 
$\reals^3$ is said to be {\em transverse} if it is everywhere transverse to the
plane field of the standard contact structure.
Parameterizing $X$ as before, a closed braid $X$ is transverse if and only if
$d\theta/dt >  - z^\prime (t)/(\rho (t))^2$ for every $t\in [0,2\pi]$. We use the symbol $TX$
to denote a transversal knot.  (For most of our discussion $TX$ will be transversal because it will also
be a braid with respect to the $z$--axis.)  Define two transverse
knots to be {\em transversally isotopic} if there is an isotopy $h\co \reals^3\times I \to
\reals^3$ such that  $h_s(X) = h(X,s) \subset \reals^3$ is transversal for every $s\in I$. 
The {\em transversal knot type} $\cTX$ of a
transverse knot $TX$ is its equivalence class under transversal isotopy.  An obvious invariant of a
transverse knot type is its {\em topological knot type} $[\cTX]_{\top}= \cX$. 

In his seminal paper \cite{Bennequin}, Bennequin showed that braids could be
useful for the study of contact structures on $S^3$ or $\reals^3$.  Recently that idea was generalized
by E Giroux \cite{Giroux}, who showed
that open book decompositions of arbitrary 3--manifolds $M^3$ might be equally important for the study of
contact structures on other 3--manifolds. Since an open book structure on a 3--manifold generalizes our braid
structures on $\reals^3$, this suggests that our work in this paper could have applications to the study of tight
contact structures on other 3--manifolds.  Bennequin proved that every transversal knot is transversally isotopic to a
closed braid, and he discovered a very important numerical invariant
$\beta(\cTX)$ (the {\em Bennequin invariant})
of transversal knot type which is natural in the braid setting. To compute it when $TX$
is also a closed braid having the $z$--axis as it axis $\axis$,
consider a regular projection of $TX$
onto the plane $z=0$. Let $e(TX)$ be the the algebraic crossing number of $TX$ at the pre-images of
the double points in the projection. Let
$b(TX)$ be the braid index of $TX$, ie, its geometric linking number with $\axis$.
Then $\beta(\cTX) = e(TX) - b(TX)$. Note that both $e(X)$ and $b(X)$ are invariants of braid isotopy, and therefore $\beta$ is well-defined on
any closed braid, ie, the closed braid need not represent a transversal knot type, and invariant under braid
isotopy.  Using topological techniques Bennequin showed that for a fixed topological knot type
$\cX$, the integer $\beta$ has an upper but not a lower bound. While $\beta(X)$ is well-defined on closed braids, it
is not an invariant of $\cX$.  Indeed,  each topological knot type $\cX$ includes
infinitely many distinct transversal knot types $\cTX$.

Are there other invariants of $\cTX$ in addition to $[\cTX]_{\top}$ and $\beta(\cTX)$?  A hint that the problem
might turn out to be quite subtle was in the paper \cite{F-T} by Fuchs and Tabachnikov, who proved that while ragbags
filled with invariants of transversal knot types $\cTX$ exist, based upon the work of
V.I. Arnold in \cite{Arnold}, they are all determined by
$[\cTX]_{\top}$ and $\beta(\cTX)$. Thus, the seemingly new invariants that many people had discovered by using
Arnold's ideas were just a fancy way of encoding $[\cTX]_{\top}$ and $\beta(\cTX)$. This lead to a definition: a transversal knot type $\cTX$ is {\em transversally simple} if it is determined by
$[\cTX]_{\top}$ and $\beta(\cTX)$. 

In \cite{B-W} purely topological techniques were used to give a purely
topological restriction on $\cX$ which would insure that any associated $\cTX$ is transversally simple.
This paper begins where the work in \cite{B-W} ended.   Our main result is to  answer the question of whether all transversal knot types are
transversally simple in the negative, by giving an explicit infinite family of examples of closed $3$--braids that define transversal knot types that are not
transversally simple. Our methods are topological and indirect and
utilize the
main result in \cite{BM-MTWS}.  After a draft of this paper was posted and submitted for publication, Etnyre and
Honda proved the existence of examples, although they could not make them explicit \cite{E-H}. Their methods, very different from ours, are
based upon techniques in contact topology. 

Here is a guide to this paper. In \fullref{subsection:block-strand diagrams and templates} we discuss our terminology and state \fullref{theorem:MTWS}, a simplified version of the MTWS.   The full statement, and the proof, are to be found in \cite{BM-MTWS}.    In \fullref{section:the MTWS and transversal isotopies} our work is directed toward the proof of \fullref{theorem:a weak transversal MTWS}, a  version of \fullref{theorem:MTWS} which holds in the special case of links that are 3--braids, in the transversal setting.   In \fullref{theorem:negative flype examples} of \fullref{section:transversal simplicity and its failure} we give our promised infinite sequence of pairs of closed 3--braids that define transversal knots types that are not transversally simple, and prove that they have the properties that we claim they have. \ms

{\bf Acknowledgments}\qua  We thank Keiko Kawamuro for her careful reading of both \cite{BM-MTWS} and this paper.   We also thank Hiroshi Matsuda for his very careful reading of an earlier version of this paper, which turned up an error.   We had claimed that a certain template had braid index 6, but he showed that it had a closed 5--braid representative.  The example, and all results that depended upon it, have been removed.

The first author acknowledges partial
support from the following sources: the US National Science
Foundation, under Grants DMS-9402988, 9705019, 9973232 and 0405586.  The second author acknowledges partial
support from the following sources: the U.S. National
Science Foundations,under grants DMS-9200881, DMS-9626884 and DMS-0306062
; and the Mathematical
Sciences Research Institute, where he was a Visiting Member during
winter/spring of 1997.

\section{The MTWS and transversal isotopies}
\label{section:the MTWS and transversal isotopies}
The goal of this section is to prove \fullref{theorem:a weak transversal MTWS}, a special and limited version of \fullref{theorem:MTWS}, the
main result in \cite{BM-MTWS}. \fullref{theorem:MTWS} is valid in the topological setting, whereas \fullref{theorem:a weak transversal MTWS} is valid in the transversal setting.   

In \fullref{subsection:block-strand diagrams and templates} our goal is to review enough of the
background so that we can state \fullref{theorem:MTWS}. Its proof is in \cite{BM-MTWS}. After that, we will
redo some parts of the proof which were given in \cite{BM-MTWS} in a transversal setting, with the goal of proving a
version of the same theorem  which applies to transversal knots.  That work begins in \fullref{subsection:the construction of transversal clasp annuli}, by redoing the `basic construction' of Section 2.1 of
\cite{BM-MTWS} in the transversal setting. This construction gives us an immersed annulus which we call a
`transversal clasp annulus'. In \fullref{subsection:braid foliations on transversal clasp annuli} we review the basic ideas about the braid foliations which
were the principle tool in \cite{BM-MTWS}, and use them to superimpose a braid foliation on the characteristic
foliation which we obtained in \fullref{subsection:the construction of transversal clasp annuli} for the
transversal clasp annulus. In \fullref{subsection:a weak transversal MTWS} we use the work in \fullref{subsection:the construction of transversal clasp annuli} and \ref{subsection:braid foliations on transversal
clasp annuli} to state and prove  \fullref{theorem:a weak transversal MTWS}.  

Our goal in writing up this part of the paper has been to make it possible for readers who are interested mainly in
the applications to be able to understand everything except parts of the proof of \fullref{theorem:a weak
transversal MTWS}.   

\subsection{Block-strand diagrams, templates and the MTWS}
\label{subsection:block-strand diagrams and templates}

\begin{figure}[ht!]
\labellist\small
\pinlabel $w$ at 126 635
\pinlabel $w$ at 235 635
\pinlabel $w$ at 371 635
\pinlabel $w$ at 475 635
\pinlabel $P$ at 134 656
\pinlabel $P$ at 243 658
\pinlabel $P$ at 379 656
\pinlabel $P$ at 483 657
\pinlabel 1 at 174 683
\pinlabel 1 at 284 683
\pinlabel 1 at 416 683
\pinlabel 1 at 524 683
\pinlabel {$\bf A$} [t] at 174 647
\pinlabel {$\bf A$} [t] at 282 649
\pinlabel {$\bf A$} [t] <1pt, 0pt> at 417 645
\pinlabel {$\bf A$} [t] <.5pt, 0pt> at 521 648
\hair 5pt
\pinlabel {positive destabilization template} [t] at 227 593
\pinlabel {negative destabilization template} [t] at 470 593
\pinlabel* {\scriptsize$+$} [l] at 153 633
\pinlabel* {\scriptsize$-$} [l] at 398 631
\endlabellist
\centerline{\includegraphics[width=.9\hsize]{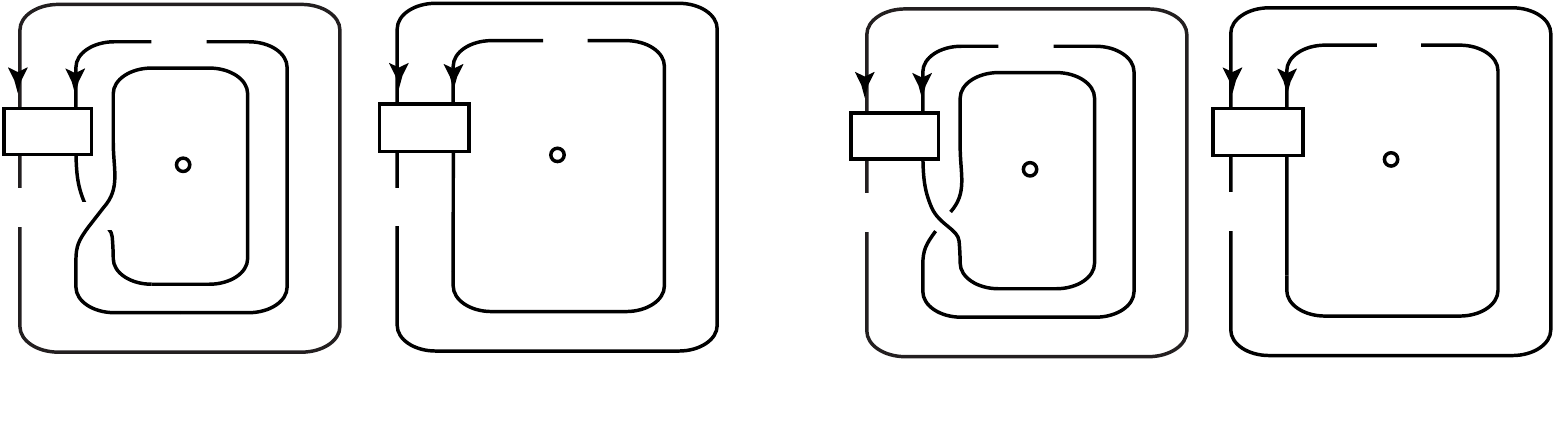}}
\caption{The two destabilization templates. Each
template is a pair of block-strand diagrams.}
\label{figure:block-strand1}
\end{figure}
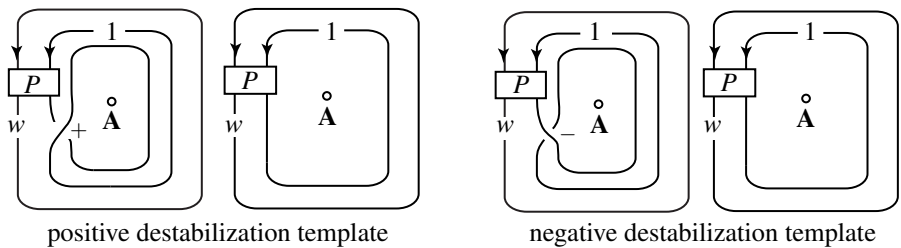

 Let $\cB$ be the collection of all braid isotopy classes of closed braid representatives of oriented
knot types in oriented $\reals^3$. Among these, consider the sub-collection
$\cB(\cX)$ of representatives of a fixed knot type $\cX$.  Among these, let $\cB_{\min}(\cX)$ be the
sub-collection of representatives whose braid index is equal to the braid index of $\cX$.  Choose any
$\Xhigh\in \cB(\cX)$  and any $\Xlow\in \cB_{\min}(\cX)$.
The classical Markov Theorem (see \cite{BM-Markov} for a proof which is
in the spirit of the work in this paper) asserts that there exists a sequence of closed braids in $\cB(\cX)$:
\begin{equation}
\label{sequence:Markov tower}
\Xhigh = X_1 \to X_2 \to \cdots \to X_r = \Xlow
\end{equation}
such that, up to braid isotopy, each $X_{i+1}$ is obtained from $X_i$ by a single stabilization or
destabilization.  
We call such a sequence a {\em Markov tower}.  The
existence of Markov towers in the transversal setting was proved by Orevkov and Shevchisin in \cite{O-S}. In the transversal setting the stabilizations and destabilizations are all required
to be positive.  That is, each $X_{i+1}$ is obtained from $X_i$ by braid isotopy and a single positive destabilization (see \fullref{figure:block-strand1}) or stabilization, where stabilization and destabilization are mutually inverse moves.   The main result in
\cite{BM-MTWS} is to find such a sequence of closed braids, like the one in \eqref{sequence:Markov tower} relating
$X_+$ to $X_-$, but not using stabilizations, so that the braid index is non-increasing.  This is done at the expense of
introducing additional moves which are non-increasing on braid index.

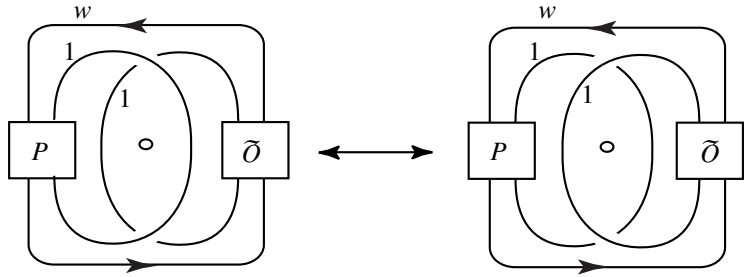
\begin{figure}[ht!]
\labellist\small
\pinlabel {\rotatebox{180}{$Q$}} at 243 643
\pinlabel {\rotatebox{180}{$Q$}} at 459 643
\pinlabel $P$ at 143 643
\pinlabel $P$ at 359 643
\pinlabel $w$ [b] at 163 702
\pinlabel $w$ [b] at 382 702
\hair 1pt
\pinlabel 1 [rb] at 161 685
\pinlabel 1 [rb] at 380 685
\pinlabel 1 [tl] at 179 673
\pinlabel 1 [tl] at 397 675
\endlabellist
\centerline{\includegraphics[scale=.8]{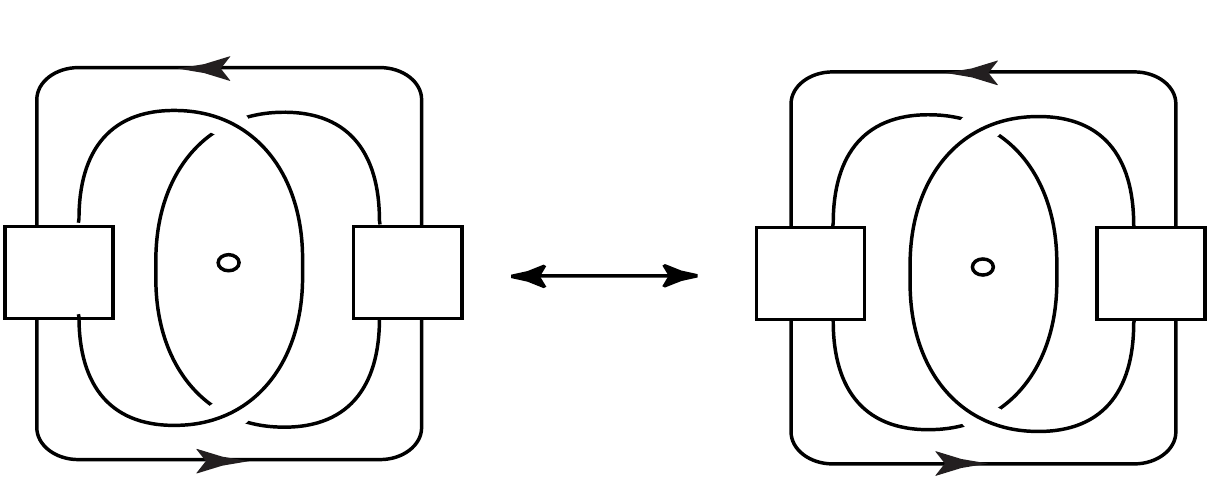}}
\caption{The exchange move template, $w\geq 2$, braid index $m\geq 4$}
\label{figure:block-strand2}
\end{figure}

\begin{figure}[ht!]
\labellist\small\hair 1pt
\pinlabel $R$ at 100 635
\pinlabel $R$ at 358 635
\pinlabel {\rotatebox{180}{$R$}} at 277 630
\pinlabel {\rotatebox{180}{$R$}} at 534 630
\pinlabel $P$ at  81 664
\pinlabel $P$ at 215 664
\pinlabel $P$ at 340 663
\pinlabel $P$ at 471 663
\pinlabel $Q$ at 82 607
\pinlabel $Q$ at 215 605
\pinlabel $Q$ at 340 606
\pinlabel $Q$ at 472 604
\pinlabel {\tiny$w$} at 92 651
\pinlabel {\tiny$w$} at 253 611
\pinlabel {\tiny$w$} at 510 610
\pinlabel {\tiny$w$} [l] at 350 649
\pinlabel {\tiny$w'$} at 124 652
\pinlabel {\tiny$w'$} at 285 607
\pinlabel {\tiny$w'$} at 383 652
\pinlabel {\tiny$w'$} at 541 609
\pinlabel {\tiny$k$} [l] at 92 622
\pinlabel {\tiny$k$} [l] at 350 621
\pinlabel {\tiny$k$} at 284 657
\pinlabel {\tiny$k$} at 541 650
\pinlabel {\tiny$k'$} at 125 607
\pinlabel {\tiny$k'$} at 253 654
\pinlabel {\tiny$k'$} at 381 606
\pinlabel {\tiny$k'$} at 510 654
\pinlabel {\tiny$+$} [l] at 150 633
\pinlabel {\tiny$+$} [r] at 227 630
\pinlabel {\tiny$-$} [l] at 407 633
\pinlabel {\tiny$-$} [r] at 481 629
\hair 5pt
\pinlabel {positive destabilization template} [t] at 188 568
\pinlabel {negative destabilization template} [t] at 445 568
\endlabellist
\centerline{\includegraphics[scale=.75]{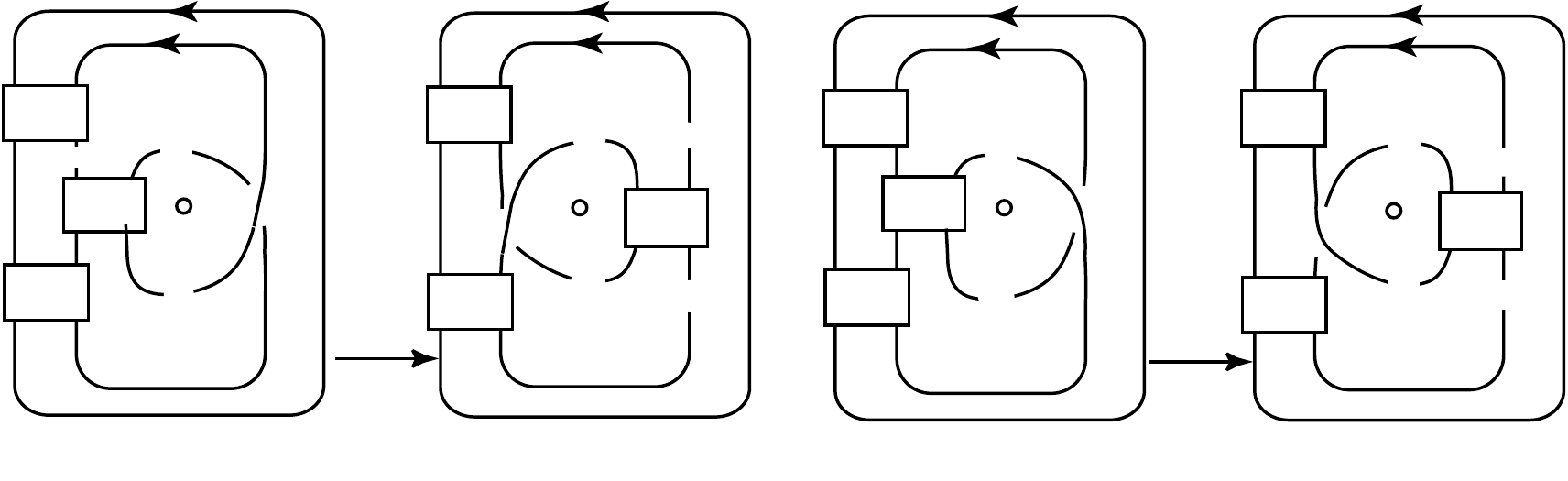}}
\caption{The two flype templates. The weights must satisfy $k'-w\geq 0$. Note that the device of weighted strands may conceal half-twists introduced
during the flype. See Figure 4 of [6].}
\label{figure:block-strand3}
\end{figure}
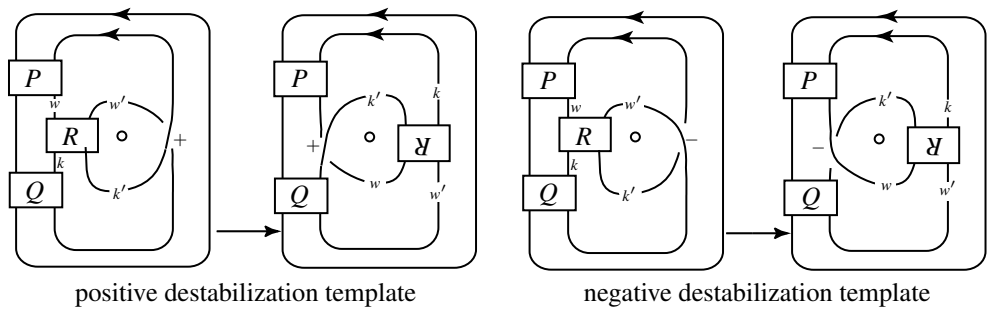

We are interested in pairs of  `block-strand diagrams' which we call `templates' (see Section 5.5,
\cite{BM-MTWS}).  The concepts are intuitive and familiar, and the reader is referred
to Figures \ref{figure:block-strand1}, \ref{figure:block-strand2} and \ref{figure:block-strand3} for
examples. 

 A {\em block} $B$ in $\reals^3-\axis$ is a 3--ball having the structure of
a $2$--disk $\Delta$ crossed with an interval $[0,1]$ such that: (i) for any fiber
$H_\theta \in \fib$ the intersection $H_\theta\cap B$ is either $\emptyset$ or $\Delta\times\{\theta\}$
for some $\theta\in [\theta_1,\theta_2]$ and, (ii) there exists some
$\theta\in [0,2\pi)$ such that $H_\theta\cap B = \emptyset$.
The disc $ \bt = B \cap H_{\theta_1} $ is the {\em top} of $B$ and the disc
$ {\bf B} = B \cap H_{\theta_2} $ is the {\em bottom} of $B$.
A {\em strand} $l$ is homeomorphic to an interval $[0,1]$ or a circle $S^1$. It is oriented 
and transverse to each fiber of $\fib$ such that its orientation agrees with the
forward direction of $\fib$.  When $l$ is homeomorphic to an interval, $\partial l = l_0 \cup l_1 $,
where $l_0$ is the {\em beginning endpoint} of $l$ and $l_1$ is the {\em ending endpoint} of $l$.

A {\em block-strand diagram}
$\cD$ is a collection of pairwise disjoint blocks 
$\{B^1, \cdots, B^{\rm k} \}$ and pairwise disjoint strands $\{l^1,\cdots,l^{\rm l}\}$
which together have the following structure:
\bi
\item[(i)] If $l^i \cap B^j \not= \emptyset$
then $l^i \cap B^j = (l_1^i \cap \bt^j ) \cup (l_0^i \cap {\bf B}^j) $ where $\bt^j$ is the top of $B^j$
and ${\bf B}^j$ is the bottom of $B^j$.  (We allow for the
possibility that either $(l_1^i \cap \bt^j )$ or $(l_0^i \cap {\bf B}^j)$ is empty.)
\item[(ii)] For each $l_0^i$ (resp.\ $l_1^i$) there is some ${\bf B}^j \subset B^j$ (resp.\ $\bt^j \subset B^j$)
such that $l_0^i \subset {\bf B}^j$ (resp.\ $l_1^i \subset \bt^j$).
\item[(iii)] For each block $B^j$ we have
$| \bt^j \cap (\cup_{1\leq i \leq {\rm l}} l_1^i ) | = | {\bf B}^j \cap (\cup_{1\leq i \leq {\rm l}} l_0^i ) | \geq 2 $.
\ei
The fact that for each $j=1,\dots,k$ there is a fiber which misses $B^j$ shows that, by rescaling, we may
find a distinguished fiber $H_{\theta_0}$ which does not
intersect any block. We define the {\em braid index}
$b(\cD)$ of the block-strand diagram
$\cD$ to be the number of times  the strands of
$\cD$ intersect the distinguished fiber $H_{\theta_0}$.  Condition (iii) above makes $b(\cD)$ well defined.

 A {\em braiding assignment} to a block-strand diagram $\cD$ 
\index{braiding assignment, block-strand diagram} 
is a choice of a braid on $m_j$ strands for each
$B^j\in\cD$. That is, we replace $B^j$ with the chosen braid, so that $B^j$ with
this braiding assignment becomes a braided tangle with $m_j$ in-strands and $m_j$
out-strands.  In this way a block strand diagram with braiding assignments gives us a closed braid
representative $X \subset \cD$
of a link type $\cX$.  We say that $X$ is {\em carried} by $\cD$.

 A {\em template} $\cT$ is a pair of block-strand diagrams
$(\Dhigh,\Dlow)$, both with blocks\break $B^1,\dots,B^{\rm k}$ and an isotopy which
takes $\Dhigh$ to $\Dlow$. That is, the blocks of $\Dhigh$ are taken onto the blocks of $\Dlow$
and the strands of $\Dhigh$ are taken onto the strands of $\Dlow$.
The diagrams $\Dhigh$ and $\Dlow$ are the {\em initial} and {\em final} block-strand
diagrams in the pair.  Thus, 
for every fixed choice of braiding assignments to the blocks
$B^1,\dots,B^{\rm k}$ the resulting closed braids $\Xhigh \subset \Dhigh$ and
$\Xlow \subset \Dlow$ represent the same oriented link type $\cX$; and the isotopy from
$\Dhigh$ to $\Dlow$ induces an isotopy from $\Xhigh$ to $\Xlow$. 

We use the  notation
$\cT(m,n)$ for the collection of all topological templates whose $\cD_+$ (resp.\ $\cD_-$) has braid index
$m$ (resp.\ $n$), where $m \geq n$.  Readers who are familiar with \cite{BM-MTWS} will note that we have changed the
definition of
$\cT(m,n)$ slightly from the definition given in Theorem 2 of \cite{BM-MTWS}, by including destabilizations, exchange
moves and flypes  in  $T(m,n)$.   We did it because it makes the statement more
concise; also we will not be concerned here, as we were in \cite{BM-MTWS}, with the actual construction of the templates,
so there is no reason to single out some templates and treat them separately from others.

 Examples of templates are the block-strand
diagram pairs which make up the templates in \fullref{figure:block-strand1},
\ref{figure:block-strand2} and \ref{figure:block-strand3}.  The strands
may carry positive integer weights, where weight $w$ means replace the strand by $w$ parallel strands in
the given projection. 

There are restrictions on the weights. The exchange move template of \fullref{figure:block-strand2} has braid
index $m = w+2$, so that it appears in every $\cT(m,n)$ with $m\geq 3.$ However, it can be replaced by braid isotopy if
$m=n=3$ (except in the special case of composite knots of braid index $3$), so we assume that $w\geq 2$. As for the
flype template of \fullref{figure:block-strand3} the restriction on the weights is needed to insure that the
flype is non-increasing on braid index. 

 A braid $X$ is said to be
{\em exchange equivalent} to $X^\prime$ if there exists an isotopy taking $X$ to $X^\prime$ that is a sequence of
braid isotopies and exchange moves.  We say $X$ is {\em exchange reducible} to $X^\prime$ if
there exists an isotopy taking $X$ to $X^\prime$ that is a sequence of braid
isotopies, exchange moves and destabilizations.

We are now ready to state the main result in \cite{BM-MTWS}, the MTWS.  We have simplified the statement, because we
do not need the full strength of the theorem proved in \cite{BM-MTWS} for present purposes.

\begin{theorem}
\label{theorem:MTWS}{{\rm (Markov's Theorem Without Stabilization \cite{BM-MTWS})}}\qua
Let $\cX$ be an oriented knot type and let $(\Xhigh , \Xlow)$ be a pair of closed braids
such that $\Xlow \in \cB_{\min}(\cX)$ and $\Xhigh \in \cB(\cX)$ with $b(\Xhigh) = m$ and $b(\Xlow) = n $.
Then there exist $\Xhigh^\prime \in \cB(\cX)$, $\Xlow^\prime \in \cB_{\min}(\cX)$,
and a finite set of templates $\cT(m,n)$ (whose membership is only dependent on $(m,n)$)
such that:
\be
\item $\Xhigh$ is exchange reducible to $\Xhigh^\prime$ and $\Xlow$ is exchange equivalent to $\Xlow^\prime$, and
\item There is a template $(\Dhigh,\Dlow)$ in the finite set $\cT(m,n)$ where, using the same braiding assignment, on $\Dhigh$
and $\Dlow$, the diagram $\Dhigh$ carries
$\Xhigh^\prime$ and the diagram $\Dlow$ carries $\Xlow^\prime$.
\ee
\end{theorem}

\subsection{The construction of transversal clasp annuli}
\label{subsection:the construction of transversal clasp annuli}

The goal in this subsection and the two that follow it  is to state and prove a modified version of \fullref{theorem:MTWS} which holds in the transversal setting. 
The proof of the MTWS in \cite{BM-MTWS} began with the construction of a  `clasp annulus'. We were given closed braids
$X_+$ and $X_-$ which represent the same topological knot type $\cX$.  The goal was to find an isotopy from
$\Xhigh\to\Xlow$ through closed braids, which would give us moves that change the closed braid, but were non-increasing on
braid index.  It was proved in Section 2.1 of
\cite{BM-MTWS}, that the trace of the isotopy from $X_+$ to $X_-$ could be assumed to 
be a {\em clasp annulus}, that is  a topological annulus
$\pa
\cong S^1
\times [0,1]$ with boundary  and an immersion $e \co  \pa \longrightarrow \ca \subset \reals^3$ such that:
\bi
\item The boundary of $\ca$  is $\Xhigh - \Xlow$. The self-intersection set of $\ca$ is the union of clasp arcs
$\{ \g^1, \cdots , \g^{\rm k} \} \subset \ca$ such that,
for each $\g^j$, we have $e^{-1}(\g^j) = \ghigh^j \cup \glow^j \subset \pa$.  Also, $\gamma_\epsilon^j$ has an 
endpoint on $e^{-1}(X_\epsilon), \ \epsilon =
\pm.$   In particular, there are no ribbon
 arcs or other more complicated self-intersections. 
\ei
We refer to $\pa$ as the
{\em pre-image annulus} (cf Section 2.1, \cite{BM-MTWS}).  The motion of pushing $\Xhigh$ across $\ca$ to
$\Xlow$ corresponds to an isotopy in $\reals^3$.  The clasp annulus $\ca$ can be positioned with respect to the
braid structure $(\fib,\axis)$ so as to have a `braid foliation' (see \fullref{subsection:braid foliations on
transversal clasp annuli} below), that encodes a Markov tower for the pair $(\Xhigh,\Xlow)$. We called this our
{\em basic construction}. 

 In this subsection we adapt the same ideas, but now in the setting of the standard contact structure and transversal
isotopies in
$\reals^3$, and construct our transversal clasp annulus.  The construction here was motivated by, but does not depend on,
that in \cite{BM-MTWS}. 

We will need to make repeated use of two results which are well-known to workers in low dimensional contact topology.  In
\cite{Bennequin}, Bennequin proved that a transversal link in the standard contact structure of $\reals^3$ could be
modified by transversal isotopy to a transversal link which is in braid position with respect
to the $z$--axis.  In fact, it is easy to see that Bennequin's original argument generalizes to the 
following statement.

\begin{lemma}[{\rm The Bennequin Trick}]
\label{lemma:the Bennequin trick}
Let $\Gamma \subset \reals^3$ be an oriented differentiable graph that is transversal to the standard
contact structure $\xi$ of $\reals^3$.  Let $e \subset \Gamma$ be an edge having a local parameterization
function $f(t)\co [0,1] \rightarrow \reals^3$ where $f(t) = (\rho(t), \theta(t), z(t))$.  If $\frac{d\theta}{dt}$ is not
positive for all $t$ then there exists a transversal isotopy of $\Gamma$ which is the identity everywhere except on a small neighborhood of
$e$ such that for a resulting parameterization $\hat{f}(t) = (\hat{\rho}(t), \hat{\theta}(t) , \hat{z}(t))$ we have
$\frac{d\hat{\theta}}{dt}$ being everywhere positive.
\end{lemma}

The Bennequin trick is analogous to the Alexander trick, which was used in Section 2.1 of
\cite{BM-MTWS}.  It essentially allows us to put a transversal graph into braid position, ie, $\Gamma \cap \axis =
\emptyset$ and 
$\frac{d\theta}{dt} > 0 $
for every edge in $\Gamma$.  The Bennequin trick will be used in conjunction with a result of Eliashberg, which
allows us to extend a transversal isotopy of links to an isotopy of the ambient space $\reals^3$, with its standard
tight contact structure.

\begin{lemma}[{\rm Eliashberg's isotopy extension lemma \cite{Eliashberg}}]
\label{lemma:contactomorphism extension}

Let $\varphi_t \co  X \rightarrow \reals^3$, $t \in [0,1]$ be an transversal isotopy of a transversal link in
the standard contact structure $\xi$.
Then there exists a contact isotopy $\alpha_t \co  (\reals^3 ,\xi) \rightarrow (\reals^3,\xi)$, $t \in [0,1]$, with
$\alpha_0 = id$ such that $\alpha_t \circ \varphi_0 = \varphi_t$ for $t \in [0,1]$.
\end{lemma}

\begin{figure}[ht!]
\labellist\hair 0pt\small
\pinlabel $TX$ [tr] at 110 546
\pinlabel $TX$ [tr] at 113 440
\pinlabel $TX$ [tr] at 093 324
\pinlabel $TX$ [tr] at 341 544
\pinlabel $TX$ [tr] at 349 441
\pinlabel $TX$ [tr] at 346 315
\pinlabel $TX_0$ [r] at 114 514
\pinlabel $TX_0$ [r] at 117 411
\pinlabel $TX_0$ [r] at 102 294
\pinlabel $TX_0$ [r] at 351 513
\pinlabel $TX_0$ [r] at 361 412
\pinlabel $TX_0$ [r] at 354 290
\pinlabel $TX_-$ [r] at 109 487
\pinlabel $TX_-$ [r] at 122 382
\pinlabel $TX_-$ [r] at 101 269
\pinlabel $TX_-$ [r] at 348 486
\pinlabel $TX_-$ [r] at 363 383
\pinlabel $TX_-$ [r] at 351 260
\hair 3pt
\pinlabel crossing [r] at 208 466
\pinlabel (a) at 208 379
\pinlabel (b) at 453 380
\pinlabel (c) at 301 245
\tiny\hair 1pt
\pinlabel $\gamma^j$ <1pt,0pt> [l] at 433 464
\pinlabel $s^j$ <0pt, 2pt> at 129 282
\pinlabel $s^{j+1}$ at 250 282
\pinlabel $s^{j'+1}$ <0.8pt,0pt> at 503 279
\pinlabel $s^{j'}$ at 379 276
\pinlabel* $\gamma^j_-$ [br] at 171 268
\pinlabel $g^i_-$ [bl] at 201 268
\pinlabel $\gamma^j_+$ at 422 320
\pinlabel $g^i_+$ at 459 320
\endlabellist
\centerline{\includegraphics[scale=.7]{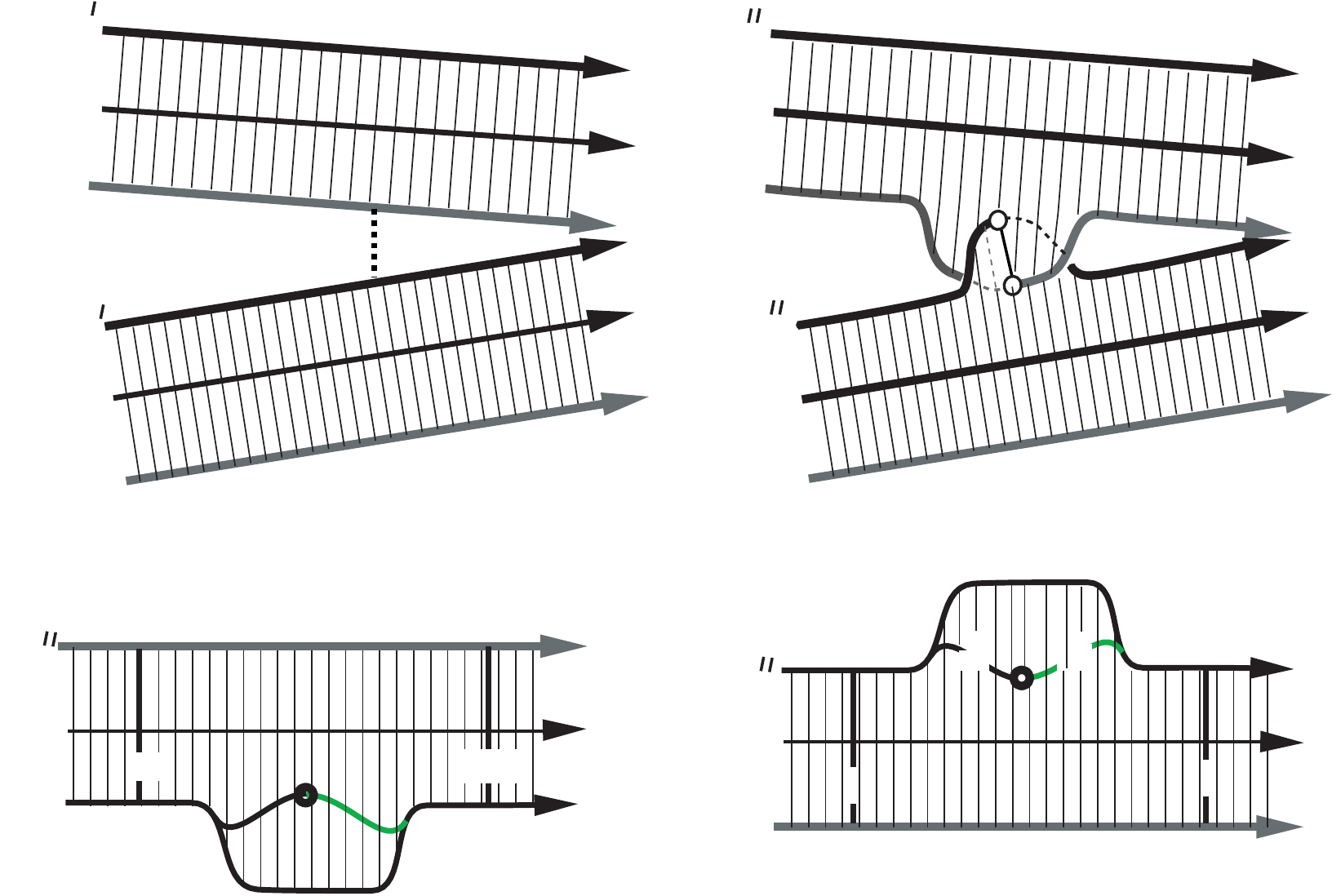}}
\caption{(a)\qua Fragments of the  embedded annulus 
$\cA^\prime$, trivially foliated, with $\partial \cA^{\prime} = TX^{\prime} - TX_-$.\qua   
(b)\qua The immersed annulus
$\cA^{\prime\prime}$, with $\partial \cA^{\prime\prime} = TX^{\prime\prime} - TX_-$.\qua 
(c)\qua The preimage of $\cA^{\prime\prime}$. (Changing the label $TX^{\prime\prime}$ to $TX_+$, sketch
(c) can also be interpreted as illustrating the preimage of $\ta$.)}
\label{figure:trans clasp ann}
\end{figure}

With the help of these two lemmas we will now develop the construction of Section 2.1 of \cite{BM-MTWS} in the contact
structure setting. \fullref{figure:trans clasp ann} may be helpful.  

Given a pair
$(\TXhigh,\TXlow)$ with
$\TXhigh, \TXlow \in \cT\cX$ and
$[\TXhigh]_{\top} , [\TXlow]_{\top} \in \cB(\cX)$, $\ta$ is an associated
{\em transversal clasp annulus} if for $\pa \cong S^1 \times [0,1]$, there exists
an immersion $e\co \pa \longrightarrow\ta \subset \reals^3$, with the oriented boundary of $\ta$ being $\TXhigh - \TXlow$, such that:
\bi
\item[(a)] The self-intersection set of the immersed annulus $\ta$ is the union of clasp arcs
$\{ \g^1, \cdots , \g^{\rm k} \} \subset \ta$ such that, for each $\g^j$, we have 
$e^{-1}(\g^j) = \ghigh^j \cup \glow^j.$  Also, $\gamma_\epsilon^j$ has an endpoint on
$e^{-1}(TX_\epsilon), \ \epsilon = \pm.$ 
\item[(b)] The characteristic foliation of $\ta$ (which comes from integrating $\xi$ in $\ta$)
is trivial.  That is, when viewing the pre-image of a leaf of the characteristic foliation in $\pa$,
it is a single arc having one endpoint on $e^{-1}(\TXhigh)$ and one endpoint on
$e^{-1}(\TXlow)$.
\item[(c)] The clasp arcs $\{ \g^1, \cdots , \g^{\rm k} \}$ are transversal arcs in the contact structure.
In particular, the arcs $\ghigh^i, \glow^i \subset \pa$ are transverse
to the pre-image of the characteristic foliation of $\ta$.  See \fullref{figure:trans clasp ann}(b) for an
example.  Each clasp arc is also in braid position, ie, transverse to each fiber $H_\theta$ of the braid structure
$\fib$ on $\reals^3$.  Here we assume that the braid axis is the $z$ axis.
\item[(d)] There exists a disjoint set of {\em extension arcs} of the clasp arcs
(\fullref{figure:trans clasp ann}(c)),
$\{g^1_\epsilon,\cdots,g^{\rm k}_\epsilon, \ \epsilon = \pm \} \subset \pa$, such that for each $i=1,\dots,k$ the union $g^i_\e \cup \gamma_\e$  is an continuous edgepath in $\pa$ that is transverse to the
trivial foliation of $\pa$ and has its two endpoints on $e^{-1}(TX_\epsilon)$. Each extension
arc is also in braid position.
\item[(e)] There exists a braid $\TXmed \subset
\{\ta \setminus (\cup_{1 \leq i \leq {\rm k}} \g^i \cup
\cup_{1 \leq i \leq {\rm k}} (g^i_+ \cup g^i_-))\}$ such that
$[\TXmed]_{\top} \in \cB(\cX)$ and $[\TXmed]_{\top}$ is a preferred longitude of both
$[\TXhigh]_{\top}$ (in the absence of $[\TXlow]_{\top}$) and $[\TXlow]_{\top}$ (in the absence of $[\TXhigh]_{\top}$).
\item[(f)] The graph 
$$\TXhigh \cup \TXlow \cup \TXmed \cup (\cup_{1 \leq i \leq {\rm k}} \g^i) \cup
(\cup_{1 \leq i \leq {\rm k}} e(g^i_+)) \cup (\cup_{1 \leq i \leq {\rm k}} e(g^i_-))$$ 
is in braid position with respect to the $z$--axis.
\item[(g)] The collection of edgepaths
$\{ (g_+^i \cup \ghigh^i),(g^i_- \cup \glow^i), \ i=1,\dots,k\} \subset \pa$ is {\em staggered}.  That is,
there exists a set of leaves
$\{s^1, \cdots , s^{\rm l} \} \subset \pa$ in the pre-image of the characteristic foliation  that are
disjoint from the $ g^i_\pm \cup \g^i_\pm $ edgepaths, $1 \leq i \leq {\rm k}$, such that
each component of $\pa \setminus (\cup_{1 \leq j \leq {\rm l}} s^j)$ is a disc containing at most a single edgepath from the set  $\{ (g^i_+ \cup \ghigh^i),(g^i_- \cup \glow^i), \ i=1,\dots,k\}$.
\ei

We now use The Bennequin Trick (\fullref{lemma:the Bennequin trick}) and the Eliashberg Lemma
(\fullref{lemma:contactomorphism extension}) to establish that for any such pair $(\TXhigh ,\TXlow)$ there always exists a transversal clasp annulus. 

\begin{proposition}
\label{proposition:transversal clasp annulus}
Given any pair of transversal closed braids $(\TXhigh,\TXlow)$ of the same transversal knot type,  there exists an associated transversal clasp annulus $\ta$ with $\partial\ta = \TXhigh - \TXlow$.
\end{proposition}

\pf We start by partitioning $\reals^3$ along the plane $z=0$, and denote the upper-half space by
$\reals^3_+$ and the lower-half space by $\reals^3_-$.  We position $\TXlow$ so that it is a
transverse braid in $\reals^3_-$.  Let $TX^\prime \in \cT\cX$ be an oriented closed braid that is a preferred
longitude of $\TXlow$.  We can visualize $TX^\prime$ being in $\reals^3_-$, in braid position, and close to $\TXlow$ with $TX^\prime - \TXlow$ the boundary of an embedded annulus
$\cA^\prime$ in the lower half-space $\reals^3_-$
whose characteristic foliation is trivial.  In \fullref{figure:trans clasp ann}(a) we see two fragments of
$\cA^\prime$. Note that $\cA'$ is oriented by the orientation on its boundary curves, so that we are looking at its negative side.   The knots $TX^\prime$ and
$\TXlow$ will be geometrically linked, except in the special case when $\cT\cX$ is the unknot.  Let
$\TXmed \subset \cA^\prime$ be a core curve of the annulus. Clearly
$\TXmed$ is also a preferred longitude of $\TXlow$ and $TX^\prime$.

Now we think of a regular projection of the link $\TXlow \sqcup TX^\prime$ onto the plane $z =0$.
Viewing this projection from the $\reals^3_+$ side, at any crossing of $\TXlow$ and $TX^\prime$,
whenever $TX^\prime$ goes under $\TXlow$ we modify it by an isotopy of $TX^\prime$, changing the
crossing. (We do not change any of the crossings between $TX^\prime$ and $\TXmed$.)  We would like to do this
`unlinking' transversally, and for that we
use a method which was used by Nancy Wrinkle in her thesis \cite{Wrinkle}. She notes:
\bi
\item [$\star$] Let
$C_\rho$ be a cylinder of radius $\rho$ in $\reals^3$ whose axis is
parallel to the $z$--axis and is
foliated by parallel Legendrian arcs of slope $-\rho^2$. Then an isotopy
of subarcs of a transverse braid
which lies on the surface of $C_\rho$, and is slid vertically along the
Legendrian leaves, is a
transverse isotopy.

\ei
So, let $p = p_i\in\reals^2$ be a double point in the projected image of
$\TXlow\cup TX'$ onto $\reals^2$,
let $p_i'$ be its preimage on $TX'$ and let
$\alpha_i$ be a neighborhood of $p_i'$ on $TX'$. For each $i$ the points $p_i$ and $p_i'$
share the same $\rho$ coordinate $\rho_i$, so they lie on the cylinder $C_{\rho_i}$. Choose a leaf
of the characteristic foliation of $C_{\rho_i}$ that passes through $p_i$ and shift the point
$p_i'$ a little bit until it also lies on the leaf. Then slide $\alpha_i$ (keeping its endpoints fixed) to
an arc $\beta_i$, thereby changing the crossing in $\TXlow\cup TX'$. Let
$TX^{\prime\prime}$ be the transverse knot
that results from repeating this procedure at each undercrossing of $TX'$,
replacing each $\alpha_i$ by $\beta_i$. 

Our unlinking isotopy will create clasp intersections
in $\cA^\prime$ as shown in \fullref{figure:trans clasp ann}.
We can assume that the clasp arcs are in braid
position. We call the resulting annulus $\cA^{\prime\prime}$. Its boundary will be $TX^{\prime\prime} - \TXlow$.
The characteristic foliation on $\cA^{\prime\prime}$ is still trivial and the image of extension arcs can be easily
chosen in $\cA^{\prime\prime}$.  Moreover, a set of leaves in the characteristic foliation can readily
chosen to demonstrate that the edgepath of clasps and their extension are staggered.  (See \fullref{figure:trans clasp ann}.)
Notice that $TX^{\prime\prime}$ and $\TXlow$ are
geometrically unlinked.  So we can modify $TX^{\prime\prime}$ by transversal isotopy so that it is contained
in $\reals^3_+$.  This transversal isotopy will also correspond to a braid isotopy in the complement of the
$z$--axis, which is the braid axis.  Finally, notice that for the pair $(X^{\prime\prime},\TXlow)$, the annulus $\cA^{\prime\prime}$ is a transversal clasp annulus.

By \fullref{lemma:the Bennequin trick} we 
know that closed oriented braids in $\reals^3$ having the $z$--axis as their braid axis are transversally isotopic to
transversal closed braids.  We consider a pair $(\TXhigh ,\TXlow)$ such that
$[\TXhigh]_{\top} \in \cB(\cX)$,
$[\TXlow]_{\top} \in \cB_{\min}(\cX)$
and $\TXhigh ,\TXlow \in \cT\cX$.
The main result in \cite{O-S} tells us that there is a {\em transversal Markov tower} for
the pair of transversal knots $(\TXhigh, \TXlow)$.  That is,
the transversal isotopy may be assumed to be via a sequence of transverse closed
braids, say:
\begin{equation}
\label{equation:Orevkov-Shevchisin}
\TXhigh = TX_1\to TX_2 \to \cdots \to TX_r = \TXlow,
\end{equation}
such that each $TX_{j+1}$ is obtained from $TX_j$ by braid isotopy and a positive stabilization or
destabilization.  Negative stabilizations and destabilizations do not occur.    

Without loss of generality,
we can assume that this isotopy is restricted to $\reals^3_+$, ie, it is the identity on
$\reals^3_-$.  By \fullref{lemma:contactomorphism extension} our transversal isotopy
via the Markov tower extends to a contact isotopy
$\Psi_0\co (\reals^3,\xi)\to (\reals^3,\xi)$, which fixes $\reals^3_-$ and, thus $\TXlow$.
This is the starting point for our transversal clasp annulus.  We initially let $\ta_0 =
\Psi_0(\cA^{\prime\prime})$. The boundary of $\ta_0$ will then be
$\Psi_0(TX^{\prime\prime} - \TXlow) = \TXhigh -\TXlow$. 
Note that since $\Psi_0$ is a contact isotopy the characteristic foliation of $\ta_0$ will be trivial,
since the foliation on $\cA^{\prime\prime}$ was trivial.
Moreover, the fact that the clasp arc and their extensions are transverse and staggered in this
trivial foliation is preserved by $\Psi_0$.
But, we are not yet done.  Although $\TXhigh$ and $\TXlow$ are in braid position, we cannot be sure that
$\TXmed$ remains in braid position under $\Psi_0$, nor do we know the positioning of the clasp arcs
$\{\g^1, \cdots , \g^{\rm k}\}$, nor the image under the immersion of the extension arcs
$\{e(g^1_+),e(g^1_-),\cdots,e(g^{\rm k}_+),e(g^{\rm k}_-)\} $.  However, if they are not in braid position, we can apply the Bennequin trick (\fullref{lemma:the Bennequin trick})
to place them in braid position via a transversal isotopy of the graph 
$$\TXhigh \cup \TXlow \cup \TXmed \cup (\cup_{1 \leq i \leq {\rm k}} \g^i) \cup
(\cup_{1 \leq i \leq {\rm k}} e(g^i_+)) \cup (\cup_{1 \leq i \leq {\rm k}} e(g^i_-))$$
which is the identity on
$\TXhigh$ and $\TXlow$.  Appealing again to \fullref{lemma:contactomorphism extension},
we extend this transversal isotopy to a contact isotopy $\Psi_1\co (\reals^3,\xi) \to (\reals^3,\xi)$.  Our needed transversal annulus associated with the pair $(\TXhigh,\TXlow)$ will be
$\ta = \Psi_1(\ta_0) = \Psi_1 \circ \Psi_0(\cA^{\prime\prime})$.
Again, the triviality of the characteristic foliation of $\cA^{\prime\prime}$ is preserved under contact isotopies.  The fact that clasp arcs and their extensions are staggered remains true under $\Psi_1$.
And, by construction $\ta$ satisfies the other conditions of a transversal clasp annulus.\endpf

\subsection{Braid foliations on transversal clasp annuli}
\label{subsection:braid foliations on transversal clasp annuli}

The principle tool in the proof of the MTWS is the study of the
{\em braid foliation} of $\ca$, the singular foliation which is determined by the intersections of
$\ca$ with the disc fibers $H_\theta$ of $\fib$.  Braid foliations have been reviewed in several places.
We refer readers who are unfamiliar with them to the review in Sections 3 and 4 of \cite{BM-MTWS}, and if
necessary thence to
\cite{B-F} for additional details. However, there are aspects of these foliations which will be needed
here and may not be so familiar, so we recall them next.  After that, we will prove \fullref{proposition:braid foliations for transversal annuli}, which describes the braid foliation on the transversal clasp annulus of the preceding section.
\begin{figure}[ht!]
\labellist
\small
\pinlabel* $X_-$ [l] at 280 665
\pinlabel* $X_-$ [l] at 529 665
\pinlabel* $X_+$ [l] at 278 547
\pinlabel* $X_+$ [l] at 529 547
\tiny
\pinlabel $+$ <.5pt,0pt> at 116 590
\pinlabel $+$ <-.3pt,0pt> at 201 623
\pinlabel $+$ <.3pt, .3pt> at 384 572
\pinlabel $+$ <.3pt, .3pt> at 457 643
\pinlabel $-$ at 201 585
\pinlabel $-$ at 116 625
\pinlabel $-$ at 385 643 
\pinlabel $-$ at 457 570
\pinlabel {type $a_-$} at 116 643
\pinlabel {type $a_+$} at 116 570
\pinlabel {type $b$} at 201 603
\pinlabel {type $s$} at 236 642
\endlabellist
\centerline{\includegraphics[scale=.7]{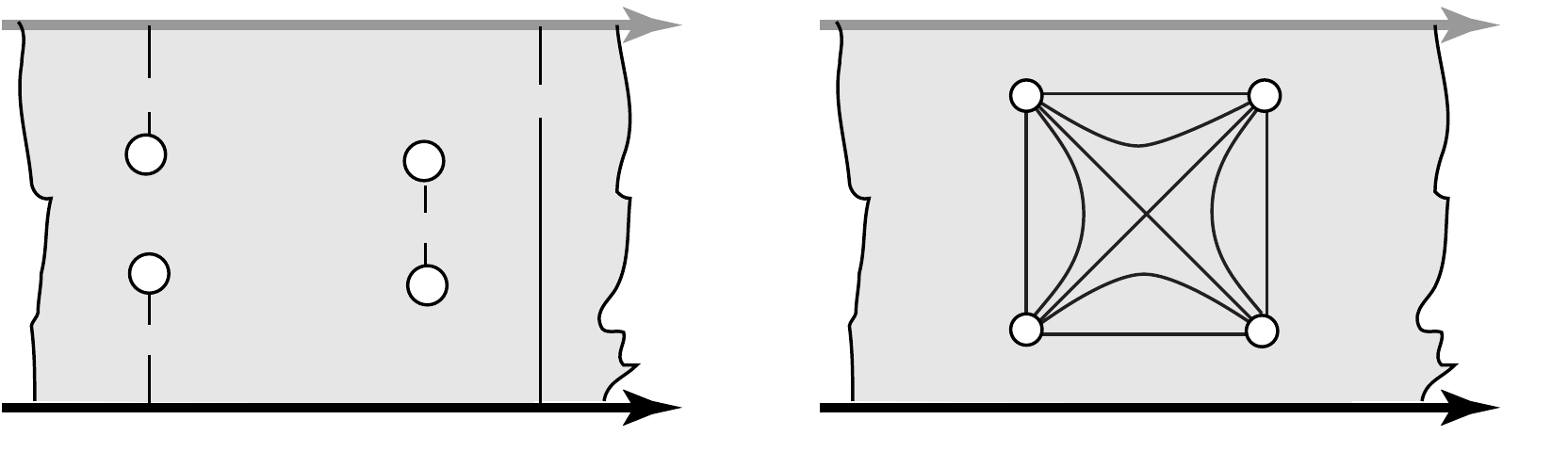}}
\caption{The left sketch illustrates the types of non-singular leaves in the foliation of $\ca$. The right sketch
illustrates a foliated neighborhood of a type $bb$--singularity.}
\label{figure:annulus1}
\end{figure}
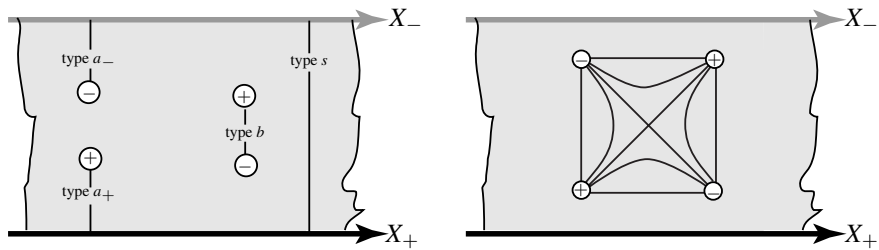

As in \cite{BM-MTWS} and \cite{B-F}: Given a braid structure
$(\fib,\axis)$ and a clasp annulus $\ca$, we arrange that $\axis$ transversally
intersects $\ca$ in finite many points.  The points of $\ca \cap \axis$ are the {\em vertices} of the
braid foliation.  Since $\ca$ has an orientation that is consistent with the orientation of $\Xhigh$, and since $\axis$
also has an orientation, each vertex of the braid foliation will have an associated parity.  We can also
assume that all but finitely many disc fibers of $H_\theta \in \fib$ will intersect $\ca$ transversally.  Of the finitely
many disc fibers that do not intersect $\ca$ transversally, each contains a single point of non-transversality that
corresponds to a saddle point.  Since $\ca$ is oriented, and since the orientations on $\axis$ and the fibers of $\fib$
are consistent,  there is also a natural parity assignment to each singularity in the braid foliation of
$\ca$. (See the left sketch in \fullref{figure:annulus1}, and see \cite{BM-MTWS} and in particular \cite{B-F} for all details). 

After an isotopy of $\ca$ with respect to the braid structure $(\fib,\axis)$ we may assume that the
non-singular leaves of the braid foliation on $\ca$ are either $s$--arcs, $a$--arcs or $b$--arcs.  
The $a$--arcs subdivide into two
groups:
$\ahigh$--arcs which have one endpoint on a positive vertex and one endpoint on $\Xhigh$; and $\alow$--arcs which have one
endpoint on a negative vertex and one endpoint on $\Xlow$. 
In the foliation the set of possible types of singularities can be listed in correspondence to the types
of leaves that are in a regular neighborhood of the singularities. 
Such a regular neighborhood is called a {\em tile}. For example, the right sketch in \fullref{figure:annulus1}
illustrates a foliated $bb$--tile.

\begin{figure}[htpb!]
\labellist
\small\hair 7pt
\pinlabel $a_+a_+$ [t] at 126 630
\pinlabel $a_+b$ [t] at 219 630
\pinlabel $a_+s$ [t] at 322 630
\pinlabel $a_-a_+$ [t] at 453 630
\pinlabel $a_-a_-$ [t] at 90 487
\pinlabel $a_-b$ [t] at 179 487
\pinlabel $a_-s$ [t] at 289 487
\pinlabel $sb$ [t] at 409 487
\pinlabel $bb$ [t] at 526 487
\tiny
\pinlabel $+$ <0pt,-.2pt> at 102 674
\pinlabel $+$ <.2pt,-.3pt> at 153 676
\pinlabel $+$ <.5pt,-.4pt> at 190 669
\pinlabel $+$ <0pt,-.4pt> at 251 669
\pinlabel $+$ at 322 632
\pinlabel $+$ at 411 656
\pinlabel $+$ <0pt,-.4pt> at 180 489
\pinlabel $+$ at 452 521
\pinlabel $+$ at 478 542
\pinlabel $+$ at 569 542
\pinlabel $-$ at 219 702
\pinlabel $-$ <0pt,-.4pt> at 65 534
\pinlabel $-$ <0pt,-.4pt> at 116 534
\pinlabel $-$ at 149 522
\pinlabel $-$ <0pt,-.4pt> at 209 522
\pinlabel $-$ <0pt,-.4pt> at 289 573
\pinlabel $-$ <0pt,-.4pt> at 526 579
\pinlabel $-$ at 526 506
\endlabellist\centerline{\includegraphics[scale=.66]{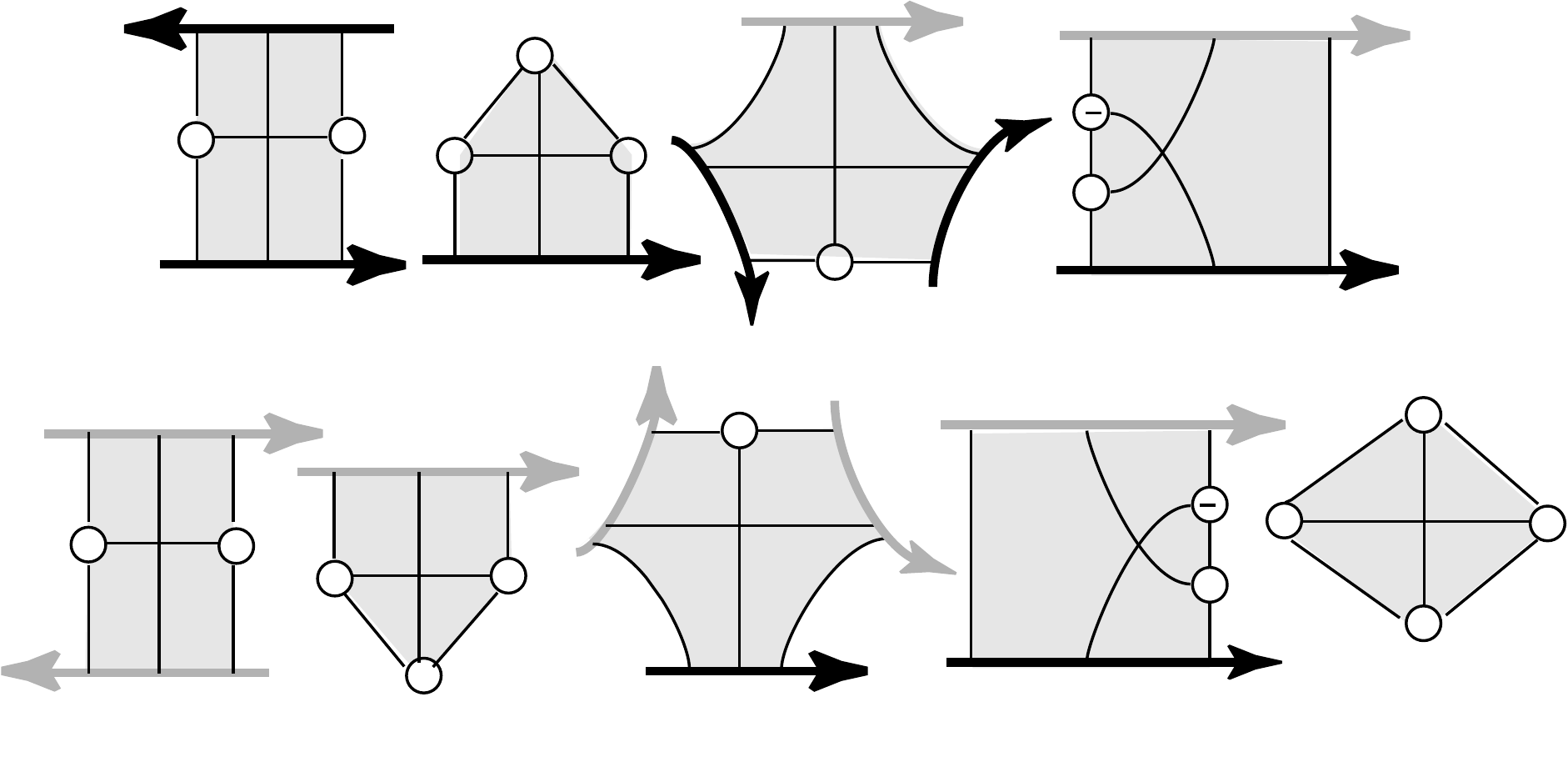}}
\caption{Possible tile types in the foliation of $\ca$.  The black (resp.\ gray) arcs
represent subarcs of $\Xhigh$ (resp.\ $\Xlow$).}
\label{figure:annulus2}
\end{figure}
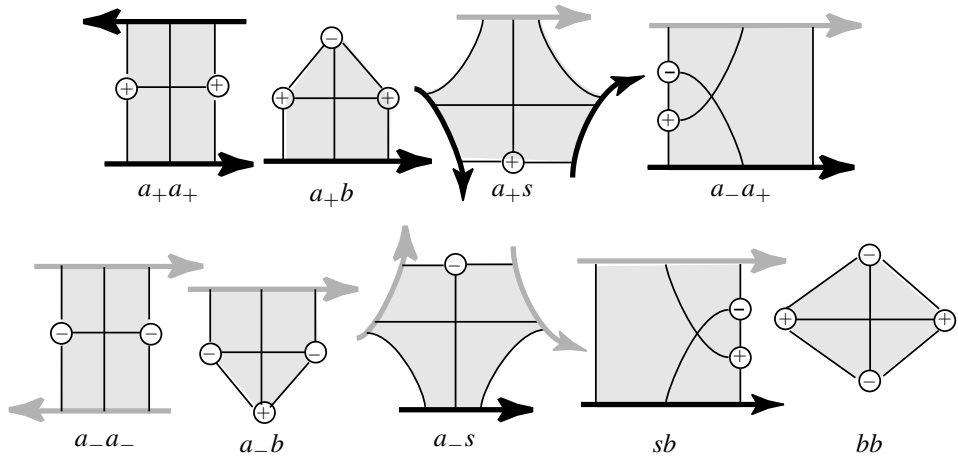

In \fullref{figure:annulus2} we show the listing of possible singularity and tile types. 

The foliation determines a cellular decomposition of $\pa$ which goes over to a cellular decomposition of $S^2$ on
shrinking the 2 boundary components to points.  Letting $V,E$ and $F$ be the number of vertices, edges and
faces, the fact that $\chi(S^2) = 2$ shows that on each component of the foliated surface $\pa$ we have
$V+2 - E + F = 2$.  Each tile has four edges and each edge is an edge of exactly 2 tiles, so that
$E = 2F$.  Combining this with the previous equation we learn that $V = F$. Let
$v_\epsilon$ (resp.\ $s_\epsilon$) be the number of vertices (resp.\ singularities) of sign $\epsilon$ in $\pa$.  
Since there is exactly one singularity
of either sign in each tile, it follows that 
\begin{equation}
\label{equation:vertices=singularities1|}
 (v_+ + v_-) = (s_+ + s_-). 
\end{equation} 
We return to the situation where $X_+$ and $X_-$ are  transversal knots $TX_+, TX_-$, also $\ta$ is a transversal clasp annulus, so that the isotopy across $\ta$ is a transversal isotopy.   We claim that the signed vertices and singularities also detect the difference between the Bennequin invariants $\beta(TX_+)$ and $\beta(TX_-)$.  For, recall that the isotopy was constructed from the transversal Markov tower of \eqref{equation:Orevkov-Shevchisin}. The moves that relate adjacent braids in that tower are braid isotopy and positive stabilizations  and destabilizations.  While braid isotopy may change the braid foliation, it does not change the numbers $v_+, v_-, s_+, s_-$, so the only issue is what happens during stabilizations and destabilizations.     Let's consider a positive destabilization.   We have a single
trivial loop around the $z$ axis, with a positive
crossing.  
We have
$d\theta > 0$ along the entire length of the loop since
we are working with a closed braid.  For a positive
crossing we have $dz
\geq 0$ throughout the loop as well.  Therefore the
inequality
$dz/d\theta > - \rho ^2$ is true for all non-zero real
values of
$\rho$.  Crossing the $z$--axis to destabilize the braid
results in at least one singular point, where
$d\theta = 0$, but if we continue to keep $dz
\geq 0$ then in the limit, as $-\rho ^2
\rightarrow 0$ from the negative real numbers,
$dz/d\theta$ goes to 
$\infty$ through the positives.  Therefore
$dz/d\theta \neq -\rho ^2$ at any stage in the
isotopy.   Therefore the isotopy is transversal.   After the isotopy the number of positive vertices will have decreased by 1 and the number of positive singularities also will have decreased by 1 so that the difference between them is zero.   Now recall that the Bennequin invariant $\beta(TX)$ of a transversal knot is the difference between its braid index and its algebraic crossing number.   Thus the Bennequin invariant is preserved during a positive destabilization.  The situation is the same for a positive stabilization.  In particular, if we study the sequence of moves that take us from $TX_+$ to $TX_-$ in \eqref{equation:Orevkov-Shevchisin} we see that:
\begin{equation}
\label{equation:vertices=singularities2}
\beta(\Xhigh) - \beta(\Xlow) = (s_+ - s_-) - (v_+ - v_-) = 0. 
\end{equation}
Using Equations \eqref{equation:vertices=singularities1|} and \eqref{equation:vertices=singularities2}, it follows that if
$\ta$ is swept out during a transversal isotopy, then: 
\begin{equation}
\label{equation:vertices=singularities}
v_+ = s_+ \ \   {\rm and} \ \ v_- = s_-.
\end{equation}
We now wish to apply our knowledge about transversal isotopies to a
specific situation. Let $X \in \cX$ be a braid in braid structure
$(\fib,\axis)$ such that $X \in \cB_{\min}(\cX)$.  Let $\TD \subset \reals^3$
be an embedded disc, with boundary $\partial(\TD)$ the union of two intervals $a_1\cup a_2$, where:
\be
\item $a_1 \subset X$.
\item $(X \setminus  a_1) \cup a_2 = X^\prime \in \cB(\cX)$.  In particular, $X \cap int(\TD) = \emptyset$ and
$X^\prime$ is isotopic to $X$ by the motion of $a_2$ to $a_1$ across $\TD$.
\item The positioning of $\TD$ in $(\fib,\axis)$ is such that the braid foliation of $\TD$ is either by parallel arcs
transverse to the boundary, or is the union of tiles from \fullref{figure:annulus2}  with the exception that we allow
the foliation of
$\TD$ to have two cusp points near the points $a_1 \cap a_2$.
\item Orienting $\TD$ so as to be consistent with the orientation of $a_2$,
the tiling of $\TD$ is such that $(v_+ - v_-) - (s_+ - s_-) = 0$.  Thus, the graph $X \cup a_2$
is in braid position and its edges have natural orientations.
\ee 
We call such a disc $\TD$ a {\em transversal disc} between $X$ and $X^\prime$. 

Before we can state our next lemma, we need one more concept about braid foliations. It is to be expected that the
union of the singular leaves will contain most of the information about the braid foliation, but in our situation
there is additional information. The parity information associated with the vertices and singularities of a braid
foliation allows us to define four graphs
$G_{+,+}, G_{+,-}, G_{-,+}$ and $G_{-,-}$.  Let $G_{\epsilon,\delta}$, where $\epsilon$ and $\delta$ are $\pm$, be
the graph that contains the
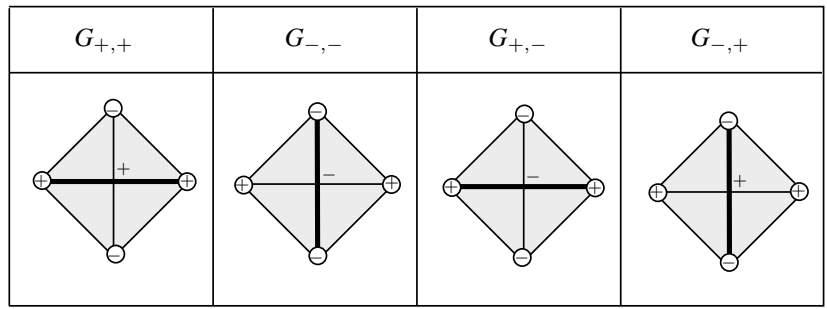
\begin{figure}[htpb]
\labellist\small\hair 7pt
\pinlabel $G_{+,+}$ [b] at 145 571
\pinlabel $G_{-,-}$ [b] at 269 571
\pinlabel $G_{+,-}$ [b] at 389 571
\pinlabel $G_{-,+}$ [b] at 508 571
\tiny\hair1pt
\pinlabel $+$ at 108 507
\pinlabel $+$ at 194 507
\pinlabel $+$ <.3pt,-.3pt> at 226 505
\pinlabel $+$ <.3pt,-.3pt> at 314 506
\pinlabel $+$ <.3pt,-.3pt> at 349 504
\pinlabel $+$ at 434 503
\pinlabel $+$ <.3pt,-.3pt> at 470 501
\pinlabel $+$ at 554 501
\pinlabel $+$ [bl] at 151 509
\pinlabel $+$ [bl] at 514 502
\pinlabel $-$ <0pt,-1pt> at 150 550
\pinlabel $-$ <0pt,-1pt> at 151 464
\pinlabel $-$ <0pt,-1pt> at 270 549
\pinlabel $-$ <0pt,-1pt> at 270 463
\pinlabel $-$ <0pt,-1pt> at 392 547
\pinlabel $-$ <0pt,-1pt> at 391 463
\pinlabel $-$ <0pt,-1pt> at 512 543
\pinlabel $-$ <0pt,-1pt> at 513 460
\pinlabel $-$ <0pt,-1pt> [bl] at 272 507
\pinlabel $-$ <0pt,-1pt> [bl] at 392 506
\endlabellist
\centerline{\includegraphics[scale=.65]{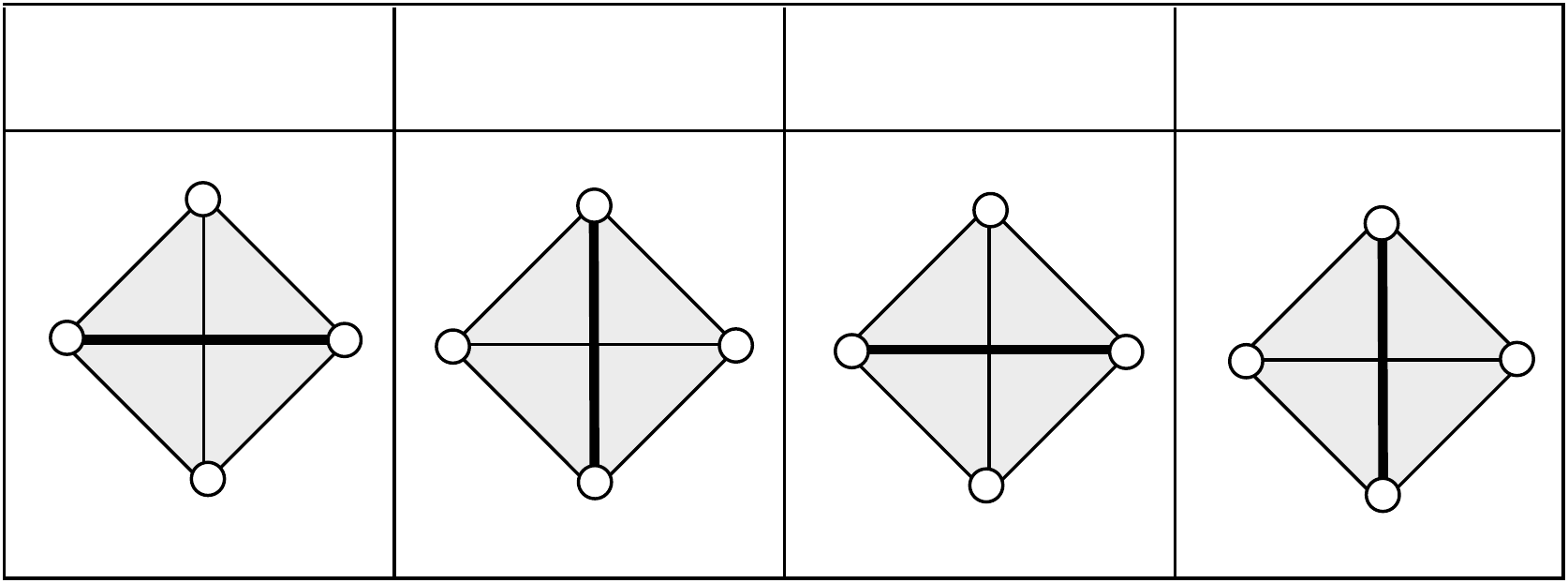}}
\caption{The graph $G_{\epsilon,\delta}$ includes all singular leaves through vertices of sign $\epsilon$
and singularities of sign $\delta$.  The thick edges illustrate the passage of
$G_{+,+}, G_{-,-}, G_{+,-}, G_{-,+}$ through a $bb$--tile.}
\label{figure:graphs1}
\end{figure}

vertices of sign $\epsilon$ and connects them together using edges that are contained in singular leaves of sign
$\delta$. By definition $G_{\epsilon,\delta}\cap G_{-\e,-\delta} =
\emptyset$. See \fullref{figure:graphs1}, which illustrates how the 4 graphs intersect a $bb$--tile. The reader should not have any difficulty in seeing how the 4 graphs pass through the tiles in \fullref{figure:annulus2}.  Observe that for all the other tiles there will be vertex endpoints of the graphs on
the boundary, because the singular leaves have endpoints on the boundary.

\begin{lemma}
\label{lemma:transversal discs}
Let $\TD$ be a transversal disc between the braids $X$ and $X^\prime$.  Then, when $X$ and $X^\prime$ are
viewed as transversal knots in the standard contact structure, the isotopy across $\TD$ is a transversal
isotopy.  In particular, if $\cC \subset G_{+,+} $ (resp.\ $G_{-,-}$) is a component of the graph in $\TD$ then
$\cC$ is a tree having a single endpoint on $a_1$ (resp.\ $a_2$).
\end{lemma}

 \pf  From \cite{B-W} we know that a positive stabilization or destabilization corresponds to a transversal
isotopy.  From Section 3.4 of \cite{BM-MTWS} we know if $\d \subset \TD$ is a regular neighborhood
of a single edge of $G_{-,-}$ which has an endpoint on $a_2$ then the isotopy of $a_2$ across
$\d$ corresponds to a positive stabilization of $a_2$.  Similarly, if $v \subset G_{+,+}$ is an
endpoint of a component and, also, a valence one vertex in the foliation of $\TD$ then from
Section 3.3 of \cite{BM-MTWS} we know we can eliminate $v$ by a positive destabilization.
Thus, if every component of $G_{-,-}$ is a tree having a single endpoint on $a_2$ we can eliminate
the graph $G_{-,-}$ through a sequence of positive stabilizations of $a_2$.  All that will remain is
the graph $G_{+,+}$.  If the components of $G_{+,+}$ are trees each having an endpoint on $a_1$ then
through a sequence of positive destabilizations (starting at vertex endpoints) we can eliminate the
components of $G_{+,+}$.  What will remain of $\TD$ will be trivially foliated and an isotopy across a
trivially foliated transverse tab corresponds to a braid isotopy in the complement of $\axis$.  But,
such braid isotopies are transversal isotopies.   Thus, to establish the main assertion of the lemma we need to
establish the assertions that the components of
our graphs $G_{\pm,\pm}$ are trees having endpoints on the appropriate $a_i$ arc of $\TD$.

The initial fact that each component of $G_{\pm,\pm}$ must be a tree follows from Theorem 3.1 of
\cite{B-F}.  We also learn from Theorem 3.1 of \cite{B-F} that if $\cC \subset G_{+,+}$ has more than one
endpoint on
$a_1$ then after a sequence of braid isotopies (which are in the form of change of fibration) and exchange moves,
we can destabilize $a_1$.  But this would imply that $X$ was not of minimal braid index.  So any component
of $G_{+,+}$ can have at most one endpoint on $a_1$.  Now notice that for those components of $G_{+,+}$ having
a single endpoint on $a_1$ the number of positive vertices equals the number of singularities.  And, for those
that don't have an endpoint on $a_1$ the number of positive vertices is always one more than the number
of singularities.  But, since $v_+ = s_+$ we have to have all components of $G_{+,+}$ having an endpoint
on $a_1$.  Finally, if $\cC \subset G_{-,-}$ is a component that has more than one endpoint on $a_2$ then there
is a path in $\cC$ that splits off a subdisc of $\TD$ that contains a component of $\cC^\prime \subset G_{+,+}$.
But, then $\cC^\prime$ could not have an endpoint on $a_1$, which contradicts what we just proved. \endpf

We need a few more definitions before we can state our main application of
\fullref{lemma:transversal discs}, ie, \fullref{proposition:braid foliations for transversal annuli}
below.  We say that a clasp arc $\g \subset \ta$ is {\em short} for $\Xhigh$ (resp.\ $\Xlow$) if there exists
a triangular disc $\DDhigh \subset \pa$ (resp.\ $\DDlow \subset \pa$)
such that $a_1 \cup a_2 \cup a_3 = \partial \DDhigh$ (resp.\ $=\partial \DDlow$)
where $a_1$ in contained in
an $a_+$--arc (resp.\ $a_-$--arc); $a_2 = \ghigh$ (resp.\ $=\glow$);
$a_3 \subset e^{-1}(\Xhigh)$ (resp.\ $e^{-1}(\Xlow)$); and, $\DDhigh$ (resp.\ $\DDlow$) is trivially foliated.

We say that a clasp arc $\glow$ (resp.\ $\ghigh$) is {\em parallel} to a component $\cC$ of
$G_{+,+}$ (resp.\ $G_{-,-}$). if there exists an edgepath $\cE \subset \cC$ and
a rectangular disc $R \subset \pa$ such that;
\be
\item $\partial R = a_1 \cup a_2 \cup a_3 \cup a_4$, where $a_1 = \glow$ (resp.\ $a_1 = \ghigh$), \ $a_2$ is a subarc
contained in a leaf of the foliation,  \ 
$a_3 \subset \Xlow$ (resp.\ $a_3 \subset \Xhigh$) and 
$a_4 = \cE$.
\item $int(R)$ does not contain any vertices or singularities of the foliation.
\ee

(A foliated neighborhood of the edgepath $\cE$ which contains $\glow$ (or $\ghigh$) and only the vertices and
singularities in
$\cE$ will play the part that tab neighborhoods of clasp arcs played in \cite{BM-MTWS}.)

\begin{proposition}
\label{proposition:braid foliations for transversal annuli}
Let $\cT\cX$ be a transversal oriented knot type and $(\TXhigh , \TXlow)$ be a pair of transversal closed braids
which represent $\cT\cX$. Assume that the braid index $b(\TXhigh) = b(\TXlow) = m$, where $m$ is minimal for all
closed braid representatives of the topological knot type $\cX$.   
Then the transversal clasp annulus $\ta$ associated with $(\TXhigh,\TXlow)$ supports a braid foliation, which
satisfies the following conditions:
\be
\item[(a)] The  leaves of the the braid foliation
are type $s$--, $\ahigh$--, $\alow$-- or $b$--arcs and the singularities correspond to those listed in \fullref{figure:annulus2}.
\item[(b)] The components of $G_{+,+}$ (resp.\ $G_{-,-}$) are trees, having a single endpoint on $\TXlow$ (resp.\ $\TXhigh$).
\item[(c)] For each clasp arc $ \g^i \subset \ta$ at least one of its associated pre-images $\ghigh^i$ or $\glow^i \subset
\pa$
is either short or is parallel to an edgepath in the pre-image of $G_{+,+}$ or $G_{-,-}$.
\ee 
\end{proposition}

\pf  Our proof will use \fullref{lemma:transversal discs}.  Let $\ta$ be a transversal clasp annulus associated with the pair
$(\TXhigh,\TXlow)$.  Inside the pre-image $\pa$ we have the pre-image of the clasp arcs
$\{\ghigh^1,\glow^1, \cdots , \ghigh^{\rm k}, \glow^{\rm k}\}$ and their extension arcs
$\{g^1_+,g^1_-, \cdots , g^{\rm k}_+, g^{\rm k}_- \}$.  Recall that we have the graph
$$\cG = \TXhigh \cup \TXlow \cup (\cup_{1 \leq i \leq {\rm k}} \g^i ) \cup ( \cup_{1 \leq i \leq {\rm k}} e(g^i_+) )
\cup ( \cup_{1 \leq i \leq {\rm k}} e(g^i_-) ), $$ and that it is both in braid position and transverse to the
standard contact structure.  By a small isotopy near the vertices of this graph we can assume that it is
differentiable.  Thus, in $\pa$ each edgepath $g^i_\pm \cup \g^i_\pm$ will split off a subdisc of $\pa$ that
has two cusps. 
\begin{figure}[htpb!]
\labellist\small
\pinlabel $TX_+$ [r] at 94 465
\pinlabel $TX_-$ [r] at 94 327
\pinlabel $TX_{os}$ [r] at 75 401
\pinlabel $\gamma^m_+$ at 133 446
\pinlabel $g^m_+$ at 176 446
\pinlabel $\gamma^i_+$ at 283 445
\pinlabel $g^i_+$ at 325 445
\pinlabel $\gamma^i_-$ at 382 347
\pinlabel $g^i_-$ <1pt,1pt> at 418 347
\pinlabel $s^1$ at 121 375
\pinlabel $s^2$ at 192 372
\pinlabel $s^3$ at 270 381
\pinlabel $s^4$ <1pt,0pt> at 335 377
\pinlabel $s^5$ <1pt,0pt> at 470 383
\endlabellist
\centerline{\includegraphics[width=.8\hsize]{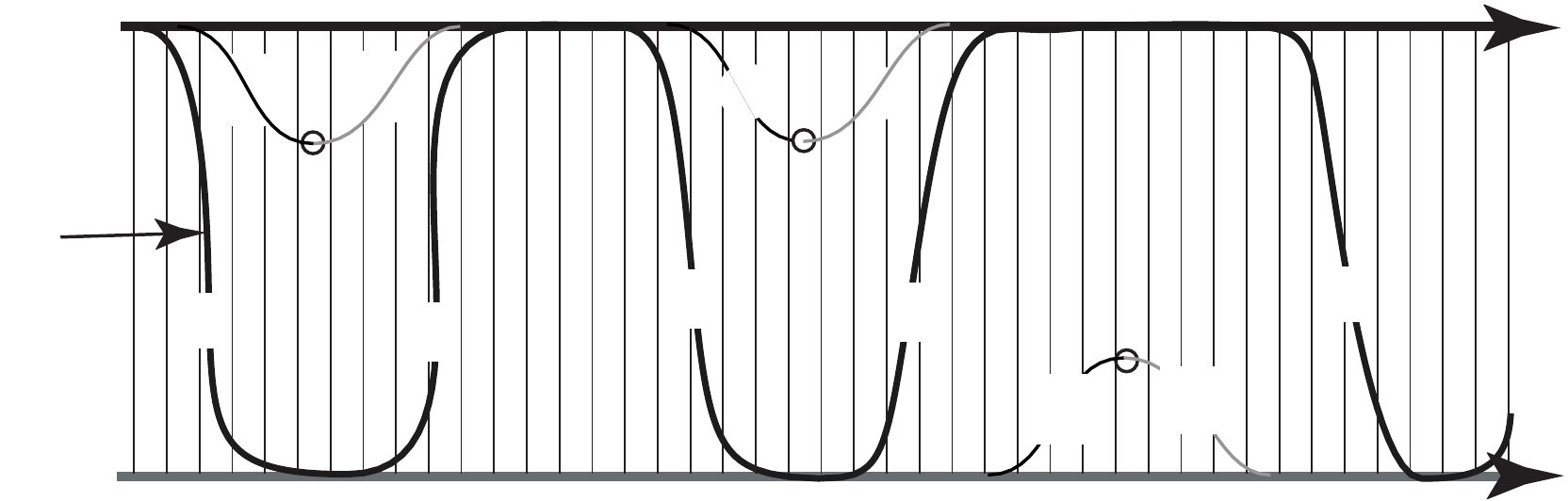}}
\caption{ The graph $\cG$, as it appears in the annulus $\pa$.}
\label{figure:oscillating}
\end{figure}
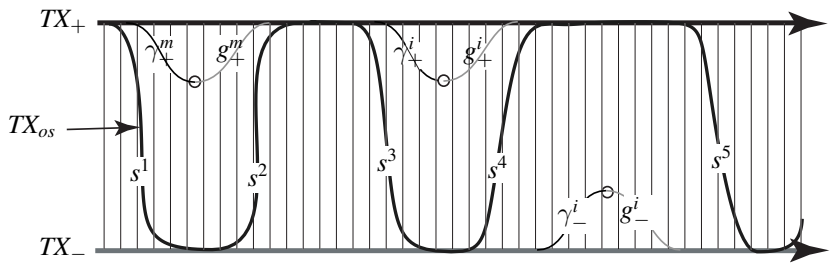

We now use the fact that the edgepaths
$\{ (g^1_+ \cup \ghigh^1),(g^1_- \cup \glow^1), \cdots , (g^{\rm k}_+ \cup \ghigh^{\rm k}),
( g^{\rm k}_- \cup \glow^{\rm k})\} \subset \pa$ are staggered.  Let $\Lambda = \{ s^1, \cdots , s^{\rm l} \} \subset
\pa$
be a set of leaves in the pre-image of the characteristic foliation of $\ta$ such that:
\be
\item[(i)] Each component of $\pa \setminus \Lambda$ is a subdisc that contains at most one $g^i_\pm \cup \g^i_\pm$
edgepath.  (We allow for the possibility that a component contains no such edgepath.) 
\item[(ii)]  By a small isotopy, the components of $\Lambda \subset \pa$ can be made transverse to the
characteristic foliation so that the graph $ \TXhigh \cup \TXlow \cup e(\Lambda) $ is: differentiable and
transverse in the standard contact structure; and, there is a smooth loop $ TX_{os} \subset \TXhigh \cup \TXlow \cup
e(\Lambda) $
that {\em oscillates} away from the clasps; ie, there exists arcs
$ \{ x^1_+, \cdots , x^{\rm l}_+ \} \subset \TXhigh $ and $\{ x^1_-, \cdots , x^{\rm l}_- \} \subset \TXlow$
such that, when the leaves of $e(\Lambda)$ are positioned transversal in $\ta$  the smooth loop
$TX_{os} = x^1_+ \cup s^1 \cup x^1_- \cup s^2 \cup x^2_+ \cup s^3 \cup x^3_- \cup \cdots \cup s^{{\rm l}}$ does
not intersect any clasp arcs.
\ee
Notice that $TX_{os} \subset \ta$ is transversally isotopic to $\TXhigh$, to $\TXmed$ and to $\TXlow$.

We now consider the expanded graph 
$$\Gamma = \TXhigh \cup \TXlow \cup TX_{os} \cup (\cup_{1 \leq i \leq {\rm k}} \g^i)
\cup ( \cup_{1 \leq i \leq {\rm k}} e(g^i_+)) \cup (\cup_{1 \leq i \leq {\rm k}} e(g^i_-)) \subset \ta$$   
The graph $\Gamma$
is a transversal graph in the standard contact structure.  Except for the portion of it that corresponds to
$e(\Lambda)$, it is also a braided graph.  We now make it entirely braided by applying first the Bennequin trick to
$e(\Lambda)$ then extending this transversal isotopy of $\Gamma$ to all of $(\reals^3 , \xi)$ by
Eliashberg's Lemma (\fullref{lemma:contactomorphism extension}).  The resulting $\ta$ will still be a transversal clasp
annulus with $TX_{os}$ now being in braid position.

We claim that the braid foliation on $\ta$ must be a tiling.  If it is not then there will exist a leaf in the
braid foliation of $c \subset \ta$ that is a circle.  But, notice that, due to the fact the $\Gamma$ is in braid
position, its edges can be consistently oriented.  Thus, $c \cap \Gamma = \emptyset$ due to orientation
considerations. We then have $c \subset \ta \setminus \Gamma$ and $c$ will therefore bound an embedded disc $ \Delta_c
\subset \ta$. Since $c \subset H_\theta$ for some $H_\theta \in \fib$ we see that it also bounds a $\Delta _\theta
\subset H_\theta$ which, without loss of generality, we can assume is inner most in $H_\theta$.  After smoothing the
corners around $c$ we will have a smooth $2$--sphere $\Delta_c \cup_c \Delta_\theta$ whose characteristic
foliation will be radial in $\Delta_\theta ( \subset \Delta_c \cup_c \Delta_\theta \cong S^2)$ with one elliptic point.  So the characteristic foliation on $\Delta_c$
must also have an elliptic point.  But, $\Delta_c \subset \ta$ and the characteristic foliation on $\ta$ is trivial.
We conclude that the circle leaf $c$ could not have existed in the braid foliation of $\ta$. Thus $\ta$ is
tiled. 

Next we establish that each clasp arc $\g^i$ is associated to the preimage of either a short clasp arc or a parallel
clasp arc in
$\pa$. So let $\cE^i_- = \glow^i \cup g^i_- \subset \pa$ be a clasp arc and its extension.  Notice that
$e(\cE^i_-) \subset \ta$ splits off an embedded disc, $\TD^i_- \subset \ta$, that has an
induced trivial characteristic foliation.  So there is a transversal
isotopy of $e(\cE^i_-)$ across $\TD^i_-$ into $\TXlow$.  Our initial goal is to establish that
$\TD^i_-$ is a transversal disc.  But, this is clear because when we consider the braid foliation of $\ta$ restricted
to $\TD^i_-$ we see that it has an induced tiling and that one of its boundary arcs is on a braid of
minimal index.  Moreover, since $\TD^i_-$ illustrates a transversal isotopy,
for the tiling of $\TD^i_-$ we have the equation $(v_+ - v_-) - (s_+ - s_-) =0$.

If the induced braid foliation on
$\TD^i_-$ is trivial, then $\glow^i$ is short.  So assume the braid foliation on $\TD^i_-$ is non-trivial.
Then we apply \fullref{lemma:transversal discs} and conclude that every component of $G_{\e,\e}, \e = \pm,$ is a
tree, also  each component of $G_{-,-}$ will have one endpoint on $e(\cE^i_-)$, and each component of $G_{+,+}$
will have one endpoint on $\TXlow$.  We can assume that the endpoints of $G_{-,-}$ on $e(\cE^i_-)$ are away
from the point $e(\glow^i \cap g^i_-) \subset \TXhigh$ which cannot be moved.  Now, by stabilizing positively
$e(\cE^i_-)$ along the components of $G_{-,-}$ in $\TD^i_-$ we can force the elimination of $G_{-,-}$ in
the resulting transversal disc, which by abuse of notation we still call $\TD^i_-$.  We now eliminate
the endpoint vertices of any component of $G_{+,+}$ in $\TD_{+,+}$ via positive destabilization, ie,
we eliminate valence one vertices of $G_{+,+}$ in $\TD^i_-$.  (See Sections 3.3 and 4.3 of \cite{BM-MTWS}.)
This can be done so long as the neighborhood of the valence one vertex does not contain the point $e(\glow^i \cap
g^i_-)$. After the elimination of all possible valence one vertices we will have the resulting $\TD^i_-$ containing
only a linear component of $G_{+,+}$ and $\glow^i$ will, thus, be parallel to an edgepath in $G_{+,+}$.

Finally, we notice that $TX_{os} \subset \ta$ decomposes $\ta$ into a union of transversal discs.
The reasoning is similar to that of why $\TD^i_-$ was a transversal disc.  Namely, any disc component
$ \TD \subset \ta \setminus TX_{os} $ is an embedded disc having
a trivial characteristic foliation coming from the characteristic foliation
of $\ta$.  Thus, an isotopy of $TX_{os}$ across $\TD$ is a transversal isotopy.  This implies that for the
braid foliation on $\ta$ restricted to $\TD$ we again have the equation $(v_+ - v_-) - (s_+ - s_-) =0$.
Applying \fullref{lemma:transversal discs} to $\TD$ we see that, since one boundary arc of
$\partial \TD$ is on a braid of minimal index, the components of $G_{\e,\e}$ are trees
having endpoints on the appropriate subarcs of the boundary of $\TD$.  This implies that the components of
$G_{\e,\e} \subset \ta$ are trees having single endpoints on $TX_{-\e}$.  Our proof is complete. \endpf

\subsection{A weak transversal MTWS}
\label{subsection:a weak transversal MTWS}

We are now ready to state and prove our transverse version of the MTWS.  

By \fullref{proposition:transversal clasp annulus} we know that there is a transversal clasp
annulus $\ta$ associated to any pair $(\TXhigh,\TXlow)$.  The characteristic foliation on $\ta$, which is trivial,
illustrates that the isotopy of $\TXhigh$ to $\TXlow$ through $\ta$ is in fact a transversal isotopy.  
We also have a second foliation on $\ta$, the singular braid foliation that is induced by the intersection of $\ta$ with
the half-planes of the standard braid fibration $\fib$ of $\reals^3$.  In \cite{BM-MTWS} this braid foliation was used
 to construct templates $(\Dhigh,\Dlow)$ of the MTWS.   Specifically, using the
graphs $G_{\epsilon,\delta}$ coming from the braid foliation, stabilizations and destabilizations along singular leaves
in the complement of clasp arcs went to make up the isotopies used in the sequences taking $\Dhigh$ to $\Dlow$.
However, our main concern, in \cite{BM-MTWS}, was in controlling the change in braid index in the passage from $\cD_+$
to $\cD_-$.  

Our concern here is somewhat different. 
Suppose we have a template $(\Dhigh,\Dlow)$ such that, when $\TXhigh$ and $\TXlow$ are viewed as braids, we see that
$\Dhigh$ carries $\TXhigh$ and $\Dlow$ carries $\TXlow$ via a common braiding assignment to the blocks of our template.
We ask whether there is a Markov tower which takes $\Dhigh$ to $\Dlow$ using only positive
stabilizations/destabilizations of braid isotopies?  In other words, by equation \eqref{equation:Orevkov-Shevchisin},
is there a transversal isotopy from $\Dhigh$ to $\Dlow$?  Having the templates in hand, we are willing to decompose
the isotopy into the moves in a markov tower.  Our only concern is whether that tower can be realized by transversal
isotopies.

Recall the notation $\cT(m,n)$ for the set of all templates in \fullref{theorem:MTWS}.  A template will 
be said to be a {\em transversal template} if there is a sequence of
positive stabilizations, positive destabilizations and braid isotopies that take $\Dhigh$ to $\Dlow$.  Our weak transversal MTWS is:

\begin{theorem}[A weak transversal MTWS]
\label{theorem:a weak transversal MTWS}
Let $\cT\cX$ be a transversal oriented knot type and let $(\TXhigh , \TXlow)$ be a pair of transversal closed braids
which represent $\cT\cX$. Assume that the braid index $b(\TXhigh) = b(\TXlow) = m$, where $m$ is minimal for all
closed braid representatives of the topological knot type $\cX$.   Then there exist
$TX_+^\prime$ and $TX_-^\prime
\in
\cB_{\min}(\cX)$, and a transversal template $(\Dhigh,\Dlow) \in \cT(m,m)$
such that:
\be
\item[(1)] $TX_+$ (resp.\ $TX_-$) is exchange equivalent to $TX'_+$ (resp.\ $TX'_-$).  In particular,
$TX_+$ (resp.\ $TX_-$) is transversally isotopic to $TX_+^\prime$ (resp.\ $TX_-^\prime$).
\item[(2)] $TX_+^\prime$ is carried by $\cD_+$ and $TX_-^\prime$ is carried by $\cD_-$, via a common braiding
assignment to $\cD_+$ and $\cD_-$.
\ee
\end{theorem}

The reader who is interested mainly in the applications to contact topology may now wish to skip ahead to the applications in \fullref{section:transversal
simplicity and its failure} in this paper.  On the other hand, the reader who has studied 
\cite{BM-MTWS} should have little difficulty in following the proof of \fullref{theorem:a weak transversal MTWS},
below.

\proof[Proof of \fullref{theorem:a weak transversal MTWS}]
Our proof will depend on the proof in
\cite{BM-MTWS} of \fullref{theorem:MTWS}, stated above. We will pinpoint the precise sections of
\cite{BM-MTWS} which are needed as they are used in our proof of \fullref{theorem:a weak transversal MTWS} by giving
section numbers and section headings (which we highlight with italics) from \cite{BM-MTWS}.  

We are given the
transversal closed braids
$(\TXhigh,\TXlow)$ which represent the transversal knot type
$\cT\cX$. By the construction in \fullref{proposition:transversal clasp annulus} we know that there is an
associated transversal clasp annulus for $(\TXhigh ,
\TXlow)$. By \fullref{proposition:braid foliations for transversal annuli}
we know that such a transversal clasp annulus will have a braid foliation such
that:  the components of $G_{+,+}$ and $G_{-,-}$ will be trees each having a single endpoint
on the appropriate boundary curve of $\ta$; and each clasp arc will either be short, or parallel to
a component of $G_{+,+}$ or $G_{-,-}$.  This is
equivalent, in the language of \cite{BM-MTWS}, to saying that the clasp arcs have been placed in tab neighborhoods. See
Section 4.3 of \cite{BM-MTWS},  {\it Construction of the tabs}.  

We ask whether the
constructions in Sections 4.4, {\it The two finger moves}, and 4.5, {\it Creating symmetric normal neighborhoods
of the clasp arcs}, of
\cite{BM-MTWS} can be used to produce symmetric normal neighborhoods for our transverse clasp arcs? It is easy to check
that the finger moves preserve the fact that our graphs are trees. Another way to see the same thing is to observe that
all the alterations to the braid foliation of $\ta$ which are used in Section 4.5 of \cite{BM-MTWS} preserve the
Bennequin number equality
$(v_+ - v_-) - (s_+ - s_-) = 0$. Thus the clasp arcs in our transversal clasp annulus can be assumed to have
symmetric normal neighborhoods, as established in Section 4.5 of \cite{BM-MTWS}.

Moving on, Section 5 of \cite{BM-MTWS} describes (in the topological setting) the ways in which one pushes 
$\Xhigh$ across $\ca$ to $\Xlow$.  These are the tools that are used to construct the templates of \cite{BM-MTWS}. We
can of course, drop the restrictions to transverse knots and consider our $\TXhigh$ and $\TXlow$ as $\Xhigh$ and
$\Xlow$, and follow the moves used in \cite{BM-MTWS} to push across $\ca$.  The question we need to address is whether
those moves can be realized by transversal isotopy? 

Clearly the complexity function of Section 5.1 of \cite{BM-MTWS} is equally valid with or
without the extra structure provided by contact topology. 

We turn our attention to Section 5.2, {\it Pushing across $\ca$ with exchange moves and destabilizations}. The
alterations that are used in the proof of Proposition 5.2.1 are (i) destabilizations of $\Xhigh$, (ii) exchange
moves and (iii) changes in foliation.  However, destabilizations of $\Xhigh$ cannot occur in our situation
because of our assumption that the topological braid index of $\TXhigh$ is minimal.  It was proved
in \cite{B-W} that exchange moves can always be realized by transversal isotopy. Note that exchange
moves on
$\TXhigh$ (resp.\ $\TXlow$) do not change  the fact that $G_{+,+}$ (resp.\ $G_{-,-}$) is a union of trees.  As for changes
in foliation, they are realized by braid isotopy. By Lemmas
\ref{lemma:the Bennequin trick} and \ref{lemma:contactomorphism extension} of this paper,  the required braid isotopies
can be realized by transversal isotopies.  So we can assume that we are in the situation at the end of Section 5.2 of
\cite{BM-MTWS}.  Putting this another way (see Section 6.1 of \cite{BM-MTWS}) we can assume that we have replaced the
original pair of closed braids $(\TXhigh, \TXlow)$ by the new pair of transversal closed braids $(\TXhigh^\prime,
\TXlow^\prime)$ of assertion (1) of \fullref{theorem:a weak transversal MTWS} of this paper.  

We move on to Section 5.3 and Section 5.4, where microflypes and flypes are used to shorten the clasp arcs and push across
$\ca$.  By the construction in \cite{BM-MTWS} all of the flypes can
be decomposed into stabilizations, microflypes, and destabilizations. The parity for the microflypes is totally
dependent on whether the clasp arcs are parallel to one of the graphs
$G_{\e,\e}$ on the one hand, or one of the graphs $G_{\e,-\e}$ on the other hand.  Since all of our clasp arcs are either
short or parallel to
$G_{+,+}$ or
$G_{-,-}$, it follows that in the braid foliation of our $\ta^\prime$ any microflype across a thin annulus will be
 realizable by a transversal isotopy.  When the microflypes are reconstituted into a larger flype,
the amalgamating conditions insure that the signs are consistent, showing that the resulting flypes are all transverse. 
The same is true in Section 5.6 of \cite{BM-MTWS}, {\it Pushing across a region with a G-flype foliation}, where
we study the amalgamation of flypes into G--flypes. The resulting G--flypes will all be transversal because our graphs
$G_{+,+}$ and $G_{-,-}$ have the needed sign data for transversal isotopy, and because without consistent signs there can
be no amalgamation of flypes into G--flypes.

We pass to G--exchange moves, studied in Section 5.7. While G--exchange moves are nothing more than inter-related
sequences of exchange moves, it seems possible that the `looping' of strands which is needed for the G--exchange moves
(illustrated nicely in the example of Figure 10 of \cite{BM-MTWS}) might not be transversal. However, when we look more
carefully at G--exchange moves, we see that they can always be replaced by sequences of stabilizations, braid isotopies
and destabilizations (see Figure 5 of \cite{BM-MTWS}), and the signs of the stabilizations and destabilizations can be
chosen to be either positive or negative for exchange moves, so again (at the expense of giving up the global nature of
the moves, which are not the issue for us here), we have transversal isotopies.   

Finally, for a motion across a standard annulus, discussed in Section 5.8 of \cite{BM-MTWS}, {\it Pushing across a
standard annulus}, we need to consider whether the components of
$G_{+,+}$ and
$G_{-,-}$ are homotopic to $S^1$ or to the unit interval. If the former, then the motion across a standard annulus will
always require a negative stabilization and destabilization. If they are homotopic to an interval, the motion will
utilize one positive stabilization, some number of exchange moves, and one positive destabilization. But since
$\ta^\prime$ satisfies the conclusions of \fullref{proposition:braid foliations for transversal annuli}, the
components of $G_{+,+}$ and $G_{-,-}$ will be homotopic to intervals and the motion will correspond to a transversal
isotopy.

The recipe of Section 5.5 of \cite{BM-MTWS}, {\it Constructing the template $(\Dhigh,\Dlow)$}, now tells us how to
construct a template $(\Dhigh,\Dlow)$ that carries  the pair $(\TXhigh^\prime,\TXlow^\prime)$.    We have proved
that  $(\Dhigh,\Dlow)$ is in fact a transversal clasp annulus, because we proved that every move that is needed can be
realized by a transversal isotopy. The proof of \fullref{theorem:a weak transversal MTWS}   is complete. \endpf

{\bf Remark}\qua  The main way in which our tranversal MTWS is restrictive and therefore `weak' is as follows. In \cite{BM-MTWS} we required
$\Xlow$ to have minimal braid index, but here we require both $\TXhigh$ and $\TXlow$ to have minimal topological braid
index. Notice that this places a very severe restriction on the transversal knot types which are covered, because
in the generic case one expects that only very special transversal knot types, for example those  which have maximal
Bennequin invariant, will be represented by closed braids with minimum topological braid index.

 \section{Transversal simplicity and its failure}
\label{section:transversal simplicity and its failure} 

 In this section we develop our examples illustrating the failure of transversal
simplicity.  In particular, we will produce examples of pairs of closed braids, $(\Xhigh,\Xlow)$ such that:
$\Xhigh,\Xlow \in \cB_{\min}(\cX)$; their associated transversal knots $\TXhigh$ and $\TXlow$ will have
$\beta(\TXhigh)=\beta(\TXlow)$; but, $\TXhigh$ is not transversally isotopic to $\TXlow$.  Our over-riding
strategy will be to employ \fullref{theorem:MTWS} and \fullref{theorem:a weak transversal MTWS} through the following line
of reasoning.

 Suppose that $\TXhigh$ and $\TXlow$ were transversally isotopic.  Then by \fullref{theorem:a weak transversal MTWS},
there must be a transversal template $(\Dhigh,\Dlow)$ that carries the pair $(\Xhigh,\Xlow)$.
By \fullref{theorem:MTWS} we know that there will only be
finitely many possible templates, and we will have judiciously chosen $(\Xhigh,\Xlow)$ so that
$(\Dhigh,\Dlow)$ is the unique template of $\cT(m,m)$ that carries $(\Xhigh,\Xlow)$.  We will
then demonstrate that $(\Dhigh,\Dlow)$ also carries links $(Y_+,Y_-)$ for which
component-wise the 
Bennequin invariant is not preserved.  This will establish that
$(\Dhigh,\Dlow)$ could not have been a transversal template.

At the end of the section we pose some open problems.

\subsection{The failure of transversal simplicity for certain closed 3--braids}
\label{subsection:the failureof transversal simplicity for certain closed 3--braids}

\begin{theorem}
\label{theorem:negative flype examples}
There exist infinitely many transversal knot types of braid index $3$
which are not transversally simple.  In particular, consider the collection of
infinitely many pairs of transversal knots $(\TXhigh,\TXlow)$ defined by the pairs of
transversal closed $3$--braids 
$$\TXhigh =  \sigma_1^{2p+1} \sigma_2^{2r} \sigma_1^{2q}\sigma_2^{-1}, \ \
\   \TXlow = \sigma_1^{2p+1} \sigma_2^{-1} \sigma_1^{2q}\sigma_2^{2r},$$ where  
$p+1\not= q\not=r$ and $p,q,r>1$
 Then, the transversal knot types
$\cTX_+$ and
$\cTX_-$ associated to each pair belong to the same topological knot type and have the same $\beta$--invariant, but
they do not represent the same transversal knot type.
\end{theorem} 

 \proof[Proof of \fullref{theorem:negative flype examples}] To begin we verify that $[{\cTX_+}]_{\top}
= [{\cTX_-}]_{\top}$ and $\beta(\TXhigh) = \beta(\TXlow)$. The closed braids 
$\TXhigh$ and
$\TXlow$ of \fullref{theorem:negative flype examples} are carried by the negative flype template
(see \fullref{figure:block-strand3}),
with the strands all assigned weight 1. The  braiding assignments to the blocks are $\sigma_1^{2p+1},
\ \sigma_2^{2r}$ and $\sigma_1^{2q}$ to $P, R$ and $Q$ respectively.  There is a template isotopy from the
left diagram to the right diagram which preserve knot or link type for any braiding assignment
to the blocks, so $[\TXhigh]_{\top}=[\TXlow]_{\top}$. If a transversal knot is defined by a
closed braid $TX$, then its $\beta$--invariant is given by the difference between the algebraic crossing
number of the diagram and the braid index, therefore $\beta(\TXhigh) = \beta(\TXlow) = 2p+2r+2q - 3$, as
claimed.

Let us assume that there is a transversal isotopy from the transverse closed braid $\TXhigh$ to
the transverse closed braid $\TXlow$. \fullref{theorem:a weak transversal MTWS} would imply that there
is a transversal
$3$--braid template that carries the braid $[\TXhigh]_{\top}$ and the braid $[\TXlow]_{\top}$.  However, we know that every
transverse isotopy is also a topological isotopy, and we know a great deal about knots that are defined by closed 3--braids, so
we can place restrictions on the topological isotopy. Let
$\cX = [\TXhigh]_{\top} = [\TXlow]_{\top}$ be the topological knot type that is defined by the transverse
closed braids $\TXhigh, \TXlow$ of \fullref{theorem:negative flype examples}.   By the main theorem in 
\cite{B-M_III}, a knot which is determined by a closed 3--braid
admits a unique braid isotopy class of closed 3--braid representatives, with the following
exceptions: 
\bi
\item the unknot, which has exactly 3 braid isotopy classes of 3--braid representatives, namely
the braid isotopy classes of the closed 3--braids
$\sigma_1^\mu\sigma_2^\tau$, where $(\mu,\tau)\in\{(1,1),(-1,-1),(1,-1)\}$; 
\item type $(2,k)$ torus knots,
where $k$ is odd and $|k|\not=1$ which have two braid isotopy classes of 3--braid representatives, namely
the conjugacy classes in
$B_3$ of $\sigma_1^k\sigma_2^\mu$, where $\mu = \pm 1$; and 
\item certain links with braid index 3 which have two
3--braid representatives, related by a positive or negative 3--braid flype. 
\ei
 From this it follows that exactly 4 templates are
needed to describe the moves of the MTWS, in the special case where
$\Xhigh$ and $\Xlow$ are closed 3--braids which define a knot of braid index at most 3: the
two destabilization templates of \fullref{figure:block-strand1} and the
two flype templates of \fullref{figure:block-strand3} with all weights 1.  In particular, every knot
or link which is represented as a closed 3--braid either has a closed braid representative which is unique
up to braid isotopy, or it is carried by one of these 4 templates. (We note that, while the exchange move
template supports knots and links of braid index 3, it can be replaced by braid isotopy for prime knots
and links of braid index at most 3.) 

Theorem 1 of \cite{Morton1} asserts that two closed 3--braids are isotopic in the complement of the braid
axis if and only if the associated elements in the 3--string braid group are conjugate in $B_3$.   Using the
solution to the conjugacy problem which is given in \cite{Murasugi} (it's due to Otto Schreier, 1936) we
selected our examples so that their conjugacy classes do not define the unknot or a type $(2,n)$--torus
knot, and therefore actually have braid index 3.  We also chose them so that we are certain that
the conjugacy classes of
$\TXhigh$ and $\TXlow$ are actually distinct.  Thus the only possibilities, in the topological setting, are
that our examples are carried by the positive or negative flype templates.  Since we know that transversal
knots that are carried by the positive flype template are transversally isotopic, we need to
choose the examples so that they are not carried by both flype templates.  For
that, we turn to \cite{K-L}, where it is shown that the closed braid
$\sigma_1^u\sigma_2^v\sigma_1^w\sigma_2^\epsilon$ admits a flype of sign $\epsilon$ and also a flype of
sign $-\epsilon$ if and only if 
$u=-\epsilon$ or $w=-\epsilon$ or $v=-2\epsilon$. So we chose our examples to avoid that possibility
too. 

We have shown that the only possibility, for the particular examples that we chose, is that the
transversal isotopy which we assumed exists from $\TXhigh$ to $\TXlow$ determines, in the topological
setting, a negative flype. Thus, our examples are uniquely carried by a single template in $\cT(3,3)$,
a negative flyping $3$--braid template.  And, by \fullref{theorem:a weak transversal MTWS} the assumption
that
$\TXhigh$ and $\TXlow$ are transversally isotopic implies this unique template is a transversal template.

We now show that this template isotopy cannot be a transversal
isotopy.  Here is the reason. A key point about the definition of a template (and this is a very strong
aspect of the MTWS) is that  for every fixed choice of braiding assignments to the blocks the resulting
closed braids represent the same oriented link type $\cX$.  For example choose the braiding assignments
$\sigma_1^3, \sigma_2^4, \sigma_1^{-5}$ to the blocks
$P,R,Q$.  This braiding assignment gives a 2--component link
$L_1\sqcup L_2$ which has two distinct isotopy classes of closed 3--braid representatives. If $L_1$ is the
component associated to the left strand entering the block $P$, then
$\beta(L_1)=-1$ and $\beta(L_2) = -3$ before the flype, but after the flype the
representative will be $\sigma_1^3 \sigma_2^{-1} \sigma_1^{-5} \sigma_2^{4}$, with 
$\beta(L_1)=-3$ and $\beta(L_2) = -1$. By Proposition
2.1.2 of \cite{Eliashberg} a  transversal isotopy of a knot/link extends to an ambient
transversal isotopy  of the 3--sphere. However, any transversal isotopy of $L_1 \sqcup L_2$
must preserve the $\beta$--invariants of the components, so no such transversal isotopy exists,
a contradiction of our assumption that $\TXhigh$ and $\TXlow$ are transversally isotopic. \endpf

\subsection{Future work and some open problems}
\label{subsection:future work and some open problems}

Here are some of the ideas and questions that occurred to us in the course of the work in this paper.

\be
\item  With the help of techniques that use only the
machinery of contact structures and characteristic foliations,  Etnyre and Honda have proved, in
\cite{E-H}, the existence of examples of pairs of transversal cable knots that are of the same knot type, have the same  Bennequin number, but are not transversally simple.   The type (2,3) cable
on a type (2,3) torus knot appears to be an example of this type,  although their methods make it difficult to establish this definitively.  
Their pair of cable knots co-bound an annulus whose characteristic
foliation has two homotopically non-trivial Legendrian closed loops.  The existence of these loops
can be seen as equivalent to the condition that the motion is across a standard annulus (Section 5.8 of
\cite{BM-MTWS}) and there are two loops in $G_{\e,\e}$.   At this writing it is an open problem to explain their non-explicit examples as explicit examples that use the machinery of the Markov Theorem Without Stabilization.    

\item Either prove or give a counterexample to the conjecture that maximal Bennequin number is always achieved at  minimum braid index.  In an earlier version of this paper we thought that we had a counterexample, but there was an error in our proof that the braids in question had braid index 6.  In fact Hirose Matsuda showed us that they had braid index 5 and not 6, so that our counterexample disappeared, and the question remains open.

\item Characterize the graphs $G_{+,+}$ and $G_{-,-}$ for an arbitrary transversal clasp annulus, that is without the
assumption that we  made here that the topological braid index is minimal.

\item Establish a transversal Markov Theorem Without Stabilization.  What should its statement be?  

\item By the work of Giroux \cite{Giroux} there is a one-to-one correspondence between
tight contact structures on arbitrary closed, orientable 3--manifolds and open book decompositions of the same
manifolds.  Noting that every open book decomposition of a closed, orientable 3--manifold can be regarded as a
branched covering space of $S^3$ with a braid structure, one expects some of the machinery in this paper and
\cite{BM-MTWS} to generalize to other 3--manifolds, yielding new connections between the topology and contact
geometry of 3--manifolds.  In this regard, the question of whether the transverse Markov theorem \cite{O-S} holds in a more general setting is probably the first question that needs to be answered.
\ee

\bibliographystyle{gtart}
\bibliography{link}

\end{document}